# A high-order Newton multigrid method with a simplified Jacobian for steady-state shallow water equations

Zhicheng Hu[a,b], Guanghan Li[a,b], Chunwu Wang[a,b], Xiaowen Wang[a,b,*]

[a] *School of Mathematics, Nanjing University of Aeronautics and Astronautics, Nanjing, 211106, China*
[b] *Key Laboratory of Mathematical Modelling and High Performance Computing of Air Vehicles (NUAA), MIIT, Nanjing, 211106, China*

**Abstract**

A high-order Newton multigrid method is proposed for steady-state shallow water flows in open channels with regular and irregular geometries. The method integrates a finite volume discretization with third-order weighted essentially non-oscillatory (WENO) reconstruction and a Newton multigrid framework with an efficient approximation of the Jacobian matrix for solving the resulting discrete system. In high-order schemes, the computational cost of Jacobian construction becomes dominant due to the wide stencil. Meanwhile, only a small fraction of the non-zero Jacobian entries exhibit large magnitudes. Based on this observation, a simplified Jacobian approximation is introduced using reduced stencils, in which selected off-stencil contributions are neglected, thereby achieving a substantial reduction in computational cost. The proposed approach is verified numerically to show significant efficiency improvement while maintaining comparable convergence behavior to that obtained with the full Jacobian approach. To further enhance performance, a geometric multigrid method incorporating a successive over-relaxation iteration as the smoother is applied to solve the linear systems arising in each Newton step. A variety of numerical experiments, including a one-dimensional smooth subcritical flow, flows over a hump, and a two-dimensional hydraulic jump over a wedge, are carried out to illustrate the third-order accuracy, efficiency, and robustness of the proposed method.

## 1. Introduction

Shallow water equations (SWEs) are widely used to model free-surface flows in rivers and coastal regions with both regular and irregular geometries [1, 2]. They provide an important mathematical framework for the analysis and prediction of flows in applications such as hydraulic engineering, environmental management, and disaster mitigation. While classical SWEs are typically formulated for idealized regular channels, extended models have been developed to account for irregular geometric effects. These variants of the SWEs are also derived from fundamental conservation laws and share a common mathematical structure for describing shallow water flows. Owing to the strong nonlinearity of the governing equations and the complicated flow features that may arise in geometrically complex settings [3, 4], the accurate and efficient numerical solution of the SWEs remains challenging, particularly for steady-state problems that require a delicate balance between flux gradients and source terms.

Numerous numerical methods have been developed for solving the SWEs, among which finite volume methods are particularly attractive due to their conservative formulation and flexibility for complex geometries. Within this framework, various

*Corresponding author.
*Email addresses:* huzhicheng@nuaa.edu.cn (Zhicheng Hu), liguanghan@nuaa.edu.cn (Guanghan Li), wangcw@nuaa.edu.cn (Chunwu Wang), wangxw@nuaa.edu.cn (Xiaowen Wang)

schemes have been proposed for shallow water flows with topography and discontinuities. For instance, a well-balanced scheme has been designed to preserve both the moving- and still-water steady states [4], while a high-resolution Godunov-type scheme has demonstrated good performance for complex two-dimensional free-surface flows [5]. To improve resolution for complicated flow structures, high-order methods based on weighted essentially non-oscillatory (WENO) reconstruction have also been widely studied [6, 7]. Comprehensive reviews of numerical methods for the SWEs can be found in [8].

Despite the considerable progress in high-order methods for the SWEs, most existing studies still compute steady-state solutions through time-dependent simulations. Although explicit time-marching schemes can eventually approach steady states, their computational efficiency is often severely restricted by the Courant-Friedrichs-Lewy (CFL) stability condition, especially for high-order discretizations. For example, the fifth-order WENO scheme [7] typically requires a CFL number below $1/12$, leading to small time-step sizes and high computational cost. Several approaches have been developed to improve robustness during time evolution, such as the multi-resolution WENO limiter in [9]. Nevertheless, these approaches do not remove the inherent efficiency limitation associated with explicit time marching toward steady states. As a result, direct steady-state solution strategies, including implicit Newton-based WENO solvers [10] and homotopy WENO methods [11], have attracted increasing attention.

The multigrid method is a widely used and highly efficient approach for steady-state computations and has demonstrated excellent convergence properties for nonlinear systems [12, 13, 14]. In particular, a Newton multigrid method (NMGM) was proposed in [15] for solving steady-state SWEs involving non-flat bottom topography. Within this framework, Newton iteration is employed to linearize the nonlinear finite volume discretization, while a geometric multigrid method with a block symmetric Gauss-Seidel smoother is used to solve the resulting linearized systems. Numerical experiments show that the method possesses good efficiency, robustness, and well-balanced properties. However, only first-order spatial accuracy was considered in [15], and relatively fine grids were therefore required to obtain accurate numerical solutions. This limitation motivates the development of high-order NMGM approaches that combine efficient multigrid acceleration with high-order spatial discretizations.

In light of these considerations, we develop a high-order NMGM for directly solving steady-state SWEs in this paper. The method combines a third-order finite volume WENO discretization with a Newton multigrid framework for solving the resulting nonlinear system. Compared with the first-order NMGM in [15], the use of high-order WENO reconstruction significantly enlarges the computational stencil and consequently increases the cost of Jacobian construction. To alleviate this difficulty, a simplified Jacobian approximation strategy based on reduced stencils is proposed for both one- and two-dimensional problems. The proposed strategy is motivated by the observation that only a small fraction of the non-zero Jacobian entries exhibit relatively large magnitudes, so that selected off-stencil contributions can be neglected without significantly affecting nonlinear convergence, thereby enhancing the efficiency of the solver. In addition, the overall efficiency is further improved by a geometric multigrid method equipped with a symmetric block successive over-relaxation (SOR) smoother. A variety of one- and two-dimensional numerical experiments are carried out to demonstrate that the proposed method achieves third-order accuracy, robustness, and substantial computational savings, while maintaining convergence behavior comparable to that of the full Jacobian formulation.

The remainder of this paper is organized as follows. In section 2, the mathematical formulations of the SWEs are introduced. Section 3 presents the high-order NMGM in the one-dimensional case, including the finite volume WENO discretization, the Newton iteration framework, the simplified Jacobian approximation strategy, and the geometric multigrid method. The extension of the method to two-dimensional problems on rectangular grids is described in section 4, followed by numerical examples in section 5. Finally, conclusions are given in the last section.

## 2. Steady-state shallow water equations

The SWEs describe the depth-averaged motion of incompressible free-surface flows under the hydrostatic pressure assumption. When a non-flat bottom topography $b(x, y)$ is considered, the two-dimensional steady-state SWEs [16] take the form

$$\begin{cases} (hu)_x + (hv)_y = 0, \\ \left(hu^2 + \dfrac{1}{2}gh^2\right)_x + (huv)_y = -ghb_x, \\ (huv)_x + \left(hv^2 + \dfrac{1}{2}gh^2\right)_y = -ghb_y, \end{cases} \tag{1}$$

where $h$ denotes the water depth, $(u, v)$ is the depth-averaged velocity vector, and $g = 9.812$ is the gravitational acceleration. The flux terms $hu^2$, $hv^2$, $huv$, and $\frac{1}{2}gh^2$ introduce nonlinearity through their quadratic dependence on the flow variables. These equations are time-independent conservation laws of mass and momentum with a source term accounting for the effect of non-flat bottom topography.

For convenience, Eqs. (1) can be rewritten into compact two-dimensional conservation form as

$$F(U)_x + G(U)_y = S(U), \tag{2}$$

with conserved variables, fluxes, and the source term given, respectively, by

$$U = \begin{pmatrix} h \\ hu \\ hv \end{pmatrix}, \quad F(U) = \begin{pmatrix} hu \\ hu^2 + \frac{1}{2}gh^2 \\ huv \end{pmatrix}, \quad G(U) = \begin{pmatrix} hv \\ huv \\ hv^2 + \frac{1}{2}gh^2 \end{pmatrix}, \quad S(U) = \begin{pmatrix} 0 \\ -ghb_x \\ -ghb_y \end{pmatrix}.$$

The system (2), though in steady state, inherits the hyperbolic character of the corresponding unsteady SWEs. For any unit vector $(n_x, n_y)^{\mathrm{T}}$ and positive $h$, the Jacobian $n_x\partial F/\partial U + n_y\partial G/\partial U$ has distinct real eigenvalues $u_n \pm \sqrt{gh}$ and $u_n$, where $u_n = un_x + vn_y$.

For flows in an open channel with bottom topography uniform in the $y$- direction, i.e. $b(x, y) \equiv b(x)$, the governing SWEs can reduce to a quasi-one-dimensional model [17], which preserves the same conservative structure as Eqs. (1) but includes additional geometric effects. Specifically, this quasi-one-dimensional model reads

$$\begin{cases} Q_x = 0, \\ \left(\dfrac{Q^2}{H} + \dfrac{1}{2}g\sigma h^2\right)_x = \dfrac{1}{2}gh^2\sigma_x - g\sigma hb_x, \end{cases} \tag{3}$$

where $\sigma(x)$ is the channel width prescribed independently of the depth, $H = h\sigma$ represents the wet cross section, and $Q = Hu$ is the mass flow rate. In the special case $\sigma(x) \equiv 1$, the above system further reduces to the classical one-dimensional SWEs, where $H = h$ and $Q = hu$. Introducing the conserved vector $U = (H, Q)^{\mathrm{T}}$, we obtain the vector form of Eqs. (3) as

$$F(U)_x = S(U), \tag{4}$$

in which $F(U) = \left(Q, \frac{Q^2}{H} + \frac{1}{2}g\sigma h^2\right)^{\mathrm{T}}$ and $S(U) = \left(0, \frac{1}{2}gh^2\sigma_x - g\sigma hb_x\right)^{\mathrm{T}}$. Here, the Jacobian $\partial F/\partial U$ admits two eigenvalues $u \pm \sqrt{gh}$, confirming the hyperbolicity of the quasi-one-dimensional system as well.

## 3. Third-order Newton multigrid method in one-dimensional case

In this section, we restrict ourselves to the one-dimensional case and introduce an efficient third-order NMGM built upon a finite volume WENO discretization.

### 3.1. Finite volume WENO discretization

Let the computational domain be divided into $N$ uniform cells of size $\Delta x$, and denote the $j$th cell by $I_j = \left[x_{j-\frac{1}{2}}, x_{j+\frac{1}{2}}\right]$, $j = 1, 2, \ldots, N$, where $x_{j\pm\frac{1}{2}}$ are the cell boundaries and $x_j = \frac{1}{2}\left(x_{j-\frac{1}{2}} + x_{j+\frac{1}{2}}\right)$ is the center. The cell average of the solution over the $j$th cell $I_j$ is defined as

$$\overline{U}_j = \frac{1}{\Delta x}\int_{I_j} U(x)\,\mathrm{d}x. \tag{5}$$

Integrating the governing one-dimensional SWEs (4) over $I_j$, we obtain the following finite volume formulation

$$\hat{F}\left(U^-_{j+\frac{1}{2}}, U^+_{j+\frac{1}{2}}\right) - \hat{F}\left(U^-_{j-\frac{1}{2}}, U^+_{j-\frac{1}{2}}\right) = \int_{I_j} S(U)\,\mathrm{d}x \triangleq S_j, \tag{6}$$

where $\hat{F}\left(U^-_{j\pm\frac{1}{2}}, U^+_{j\pm\frac{1}{2}}\right)$ are the numerical fluxes across the cell interfaces $x_{j\pm\frac{1}{2}}$, with $U^-_{j\pm\frac{1}{2}}$ and $U^+_{j\pm\frac{1}{2}}$ denoting the reconstructed left and right states, respectively. Various approximate Riemann solvers, such as the local Lax-Friedrichs (LLF) flux [18], and the Harten-Lax-van Leer (HLL) flux [4], can be employed to evaluate them. The specific choice adopted in our numerical experiments will be given in section 5.

To achieve high-order accuracy, well-resolved interface values $U^\pm_{j\pm\frac{1}{2}}$ are required. In this work, we employ a third-order WENO reconstruction [19, 20] to obtain them. For example, according to the WENO procedure, the stencil $\left\{I_{j-1}, I_j, I_{j+1}\right\}$, consisting of two sub-stencils $T_0 = \left\{I_j, I_{j+1}\right\}$ and $T_1 = \left\{I_{j-1}, I_j\right\}$ as illustrated in Figure 1a, is used to reconstruct $U^+_{j-\frac{1}{2}}$ and $U^-_{j+\frac{1}{2}}$. Specifically, these interface values are expressed as convex combinations of sub-stencil reconstructions, that is,

$$U^+_{j-\frac{1}{2}} = \sum_{r=0}^{1} \omega_r U^{(r)}_{j-\frac{1}{2}}, \quad U^-_{j+\frac{1}{2}} = \sum_{r=0}^{1} \tilde{\omega}_r U^{(r)}_{j+\frac{1}{2}}. \tag{7}$$

Here, $U^{(r)}_{j\pm\frac{1}{2}}$ are the candidate linear reconstructions obtained from the sub-stencil $T_r$ by

$$\begin{aligned} U^{(0)}_{j-\frac{1}{2}} &= \frac{3}{2}\overline{U}_j - \frac{1}{2}\overline{U}_{j+1}, \quad U^{(1)}_{j-\frac{1}{2}} = \frac{1}{2}\overline{U}_{j-1} + \frac{1}{2}\overline{U}_j, \\ U^{(0)}_{j+\frac{1}{2}} &= \frac{1}{2}\overline{U}_j + \frac{1}{2}\overline{U}_{j+1}, \quad U^{(1)}_{j+\frac{1}{2}} = -\frac{1}{2}\overline{U}_{j-1} + \frac{3}{2}\overline{U}_j. \end{aligned}$$

The corresponding coefficients $\omega_r$ and $\tilde{\omega}_r$ are normalized nonlinear weights defined as

$$\omega_r = \frac{\alpha_r}{\alpha_0 + \alpha_1}, \quad \tilde{\omega}_r = \frac{\tilde{\alpha}_r}{\tilde{\alpha}_0 + \tilde{\alpha}_1}, \quad r = 0, 1,$$

where

$$\alpha_r = \frac{d_r}{(\zeta + \beta_r)^2}, \quad \tilde{\alpha}_r = \frac{\tilde{d}_r}{(\zeta + \tilde{\beta}_r)^2}, \tag{8}$$

in which $d_r, \tilde{d}_r$ are the linear weights ensuring optimal accuracy in smooth regions, and $\beta_r, \tilde{\beta}_r$ are the smoothness indicators. On a uniform grid, they are given by

$$\begin{aligned} &d_0 = \tilde{d}_1 = \frac{2}{3}, \quad d_1 = \tilde{d}_0 = \frac{1}{3}, \\ &\beta_0 = \tilde{\beta}_0 = (\overline{U}_{j+1} - \overline{U}_j)^2, \quad \beta_1 = \tilde{\beta}_1 = (\overline{U}_j - \overline{U}_{j-1})^2. \end{aligned}$$

In addition, $\zeta = 10^{-6}$ is a small parameter introduced to prevent division by zero.

A high-order approximation of the source term is also essential to ensure the overall third-order accuracy of the discretization. To this end, the two-point Gauss-Legendre quadrature is applied to the integral on the right-hand side of (6), yielding

$$S_j = \int_{I_j} S(U)\,\mathrm{d}x \approx \int_{I_j} S(U_h)(x)\,\mathrm{d}x \approx \frac{\Delta x}{2}\left[S(U_h)\left(x_j - \frac{\Delta x}{2\sqrt{3}}\right) + S(U_h)\left(x_j + \frac{\Delta x}{2\sqrt{3}}\right)\right], \tag{9}$$

(a) Sub-stencils $T_0$ and $T_1$ for reconstruction of $U^+_{j-\frac{1}{2}}$ and $U^-_{j+\frac{1}{2}}$.

(b) The required stencil for the third-order discretization over $I_j$.

Figure 1: Stencil configurations for one-dimensional discretization.

where $U_h(x)$ is the cubic polynomial interpolated from the reconstructed interface values $U^+_{j-\frac{3}{2}}$, $U^+_{j-\frac{1}{2}}$, $U^-_{j+\frac{1}{2}}$, and $U^-_{j+\frac{3}{2}}$.

Collecting all the above components, we can readily see that the present third-order discretization (6) over $I_j$ involves five cell averages from the stencil illustrated in Figure 1b. For the sake of clarity, we prefer to rewrite (6) in a residual form as

$$R_j\left(\overline{\boldsymbol{U}}_j\right) \triangleq \hat{F}\left(U^-_{j+\frac{1}{2}},\, U^+_{j+\frac{1}{2}}\right) - \hat{F}\left(U^-_{j-\frac{1}{2}},\, U^+_{j-\frac{1}{2}}\right) - S_j = 0, \tag{10}$$

where $R_j$ denotes the local residual over $I_j$, and $\overline{\boldsymbol{U}}_j = \left(\overline{U}_{j-2}, \overline{U}_{j-1}, \overline{U}_j, \overline{U}_{j+1}, \overline{U}_{j+2}\right)$ collects all the cell averages on which $R_j$ depends.

### 3.2. Regularized Newton iteration

The discrete system, assembled from the residual equations (10) for all cells, is nonlinear with respect to the unknown variable $\overline{\boldsymbol{U}}$, which collects all cell averages $\overline{U}_j$. To solve this nonlinear system, Newton's method is considered here. It follows that the residual equation (10) over $I_j$ is linearized around the current approximation $\overline{\boldsymbol{U}}^{(n)}$ as

$$\sum_{i=j-2}^{j+2} \left(\frac{\partial R_j}{\partial \overline{U}_i}\right)^{(n)} \delta U_i^{(n)} = -R_j^{(n)}, \tag{11}$$

where $\delta U_i^{(n)}$ denotes the increment to $\overline{U}_i$ at the $n$th iteration, and $R_j^{(n)} = R_j\left(\overline{\boldsymbol{U}}_j^{(n)}\right)$ is the corresponding local residual. Due to the nonlinear flux and the WENO reconstruction procedure, the analytical evaluation of the Jacobian entries $\left(\partial R_j / \partial \overline{U}_i\right)^{(n)}$ is intractable. To circumvent this analytical difficulty and keep the formulation independent of the specific numerical flux, these Jacobian entries are computed by numerical differentiation. In particular, each $\left(\partial R_j / \partial \overline{U}_i\right)^{(n)}$ is a $2 \times 2$ matrix, and its $(l, m)$th element is approximated as

$$\left(\frac{\partial R_j}{\partial \overline{U}_i}\right)^{(n)}_{lm} \approx \frac{R_j^l\left(\overline{\boldsymbol{U}}_j^{(n)} + \epsilon \boldsymbol{e}_{im}\right) - R_j^l\left(\overline{\boldsymbol{U}}_j^{(n)}\right)}{\epsilon}, \quad l, m = 1, 2, \tag{12}$$

where $R_j^l$ denotes the $l$th component of the residual vector, $\epsilon$ is a small perturbation, and $\boldsymbol{e}_{im} = \left(\delta_{j-2,i}, \delta_{j-1,i}, \delta_{ji}, \delta_{j+1,i}, \delta_{j+2,i}\right) e_m$, in which $\delta_{ki}$ is the Kronecker delta symbol and $e_m$ is the unit vector with the $m$th component equal to 1 and 0 elsewhere.

For the hyperbolic problem, the global linear system assembled from (11) for all $j$ can be singular and thus requires regularization. Following the strategy adopted in [12, 13, 21], the $\ell^1$-norm of the local residual, $\left\|R_j^{(n)}\right\|_{\ell^1}$, is employed to regularize (11), leading to

$$\alpha \left\|R_j^{(n)}\right\|_{\ell^1} \mathbf{I}\, \delta U_j^{(n)} + \sum_{i=j-2}^{j+2} \left(\frac{\partial R_j}{\partial \overline{U}_i}\right)^{(n)} \delta U_i^{(n)} = -R_j^{(n)}, \tag{13}$$

where $\alpha$ is a positive regularization parameter and $\mathbf{I}$ is the identity matrix.

Once all increments $\delta U_j^{(n)}$ are obtained by solving the above regularized linear system, the cell averages are updated by

$$\overline{U}_j^{(n+1)} = \overline{U}_j^{(n)} + \tau\, \delta U_j^{(n)}, \quad j = 1, \ldots, N, \tag{14}$$

where $\tau \in (0, 1]$ is a damping parameter. This completes one iteration of Newton's method. Repeating the iteration until convergence constitutes a regularized Newton-based solver, which is summarized in **Algorithm 1**. In the following subsections, the efficient construction of an approximate Jacobian matrix and the solution strategy for the resulting linear system will be discussed in detail.

**Algorithm 1** Regularized Newton-based solver

1: Initialize $\overline{U}_j^{(0)}$, $j = 1, \ldots, N$; specify the tolerance *Tol* and the maximum number of iterations $n_{\max}$; set $n = 0$;
2: **repeat**
3: Compute the Jacobian entries via numerical differentiation (12);
4: Regularize the linear system (11) using the $\ell^1$-norm of the local residual $R_j^{(n)}$;
5: Solve the regularized system resulting from (13) to obtain the increment $\delta U_j^{(n)}$ on each cell;
6: Update the cell averages $\overline{U}_j^{(n+1)}$ according to (14);
7: Evaluate the global residual $res = \frac{1}{N} \sum_{j=1}^{N} \left\| R_j^{(n+1)} \right\|_{\ell^1}$;
8: Set $n \leftarrow n + 1$;
9: **until** $res < Tol$ or $n \geq n_{\max}$.

### 3.3. *Efficient approximation of the Jacobian matrix*

The full Jacobian matrix, with its $(j, i)$th block given by $\alpha \left\| R_j^{(n)} \right\|_{\ell^1} \delta_{ji} \mathbf{I} + \left( \partial R_j / \partial \overline{U}_i \right)^{(n)}$, is a block-pentadiagonal matrix of size $N \times N$. It is hereafter denoted by $\boldsymbol{J}5 = \left( J_{ji} \right)_{N \times N}$. Clearly, its evaluation accounts for a major portion of the total computational cost in each Newton iteration.

To improve computational efficiency without compromising the convergence of the Newton iteration, the magnitudes of all nonzero blocks in $\boldsymbol{J}5$, measured by the entrywise $\ell^1$-norm, are first examined. Numerical observations show that, in most situations, the outermost diagonal blocks are of considerably smaller magnitude than those along the central diagonals. A typical distribution of these block magnitudes for $N = 10$ is illustrated in Figure 2a. More quantitatively, for $\boldsymbol{J}5$ of size $96 \times 96$, obtained from the first numerical experiment presented in section 5, the magnitudes of the central tridiagonal blocks consistently exceed 1, whereas those of the outermost diagonal blocks fall below 1. Notably, nearly 69% of the outermost diagonal blocks exhibit magnitudes even smaller than $1 \times 10^{-4}$.

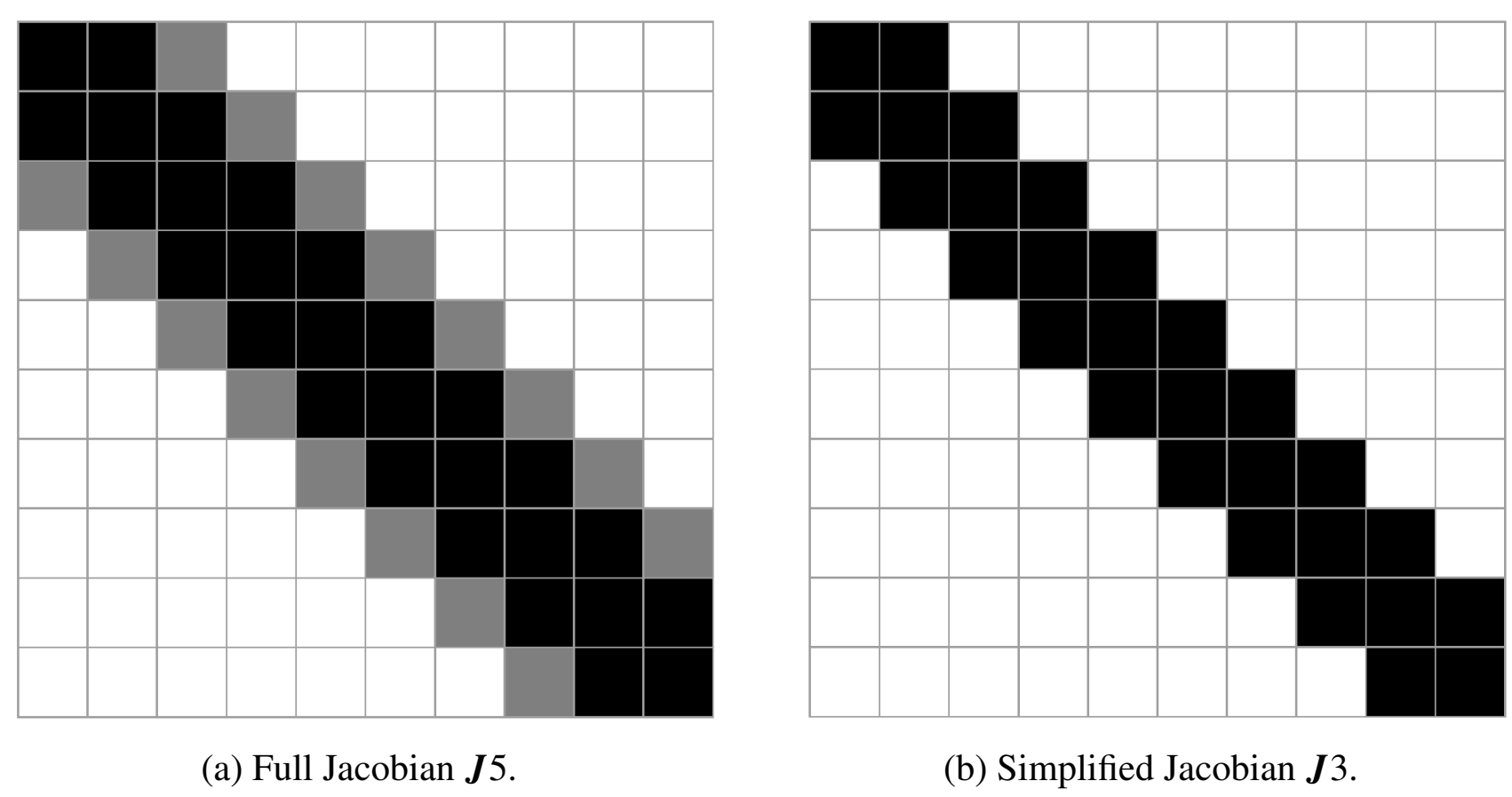

(a) Full Jacobian $\boldsymbol{J}5$.   (b) Simplified Jacobian $\boldsymbol{J}3$.

Figure 2: Typical patterns of block magnitudes in the Jacobian matrix for $N = 10$. Black squares indicate blocks with magnitude $> 1$ and gray ones $< 1$.

Based on the above observations, it is reasonable to expect that the outermost diagonal blocks contribute very little to the overall Jacobian and can thus be neglected in its construction. This directly leads to a simplified block-tridiagonal Jacobian, denoted by $\boldsymbol{J}3$, with the pattern of block magnitudes illustrated in Figure 2b. Obviously, such a simplification reduces the computational cost of evaluating the Jacobian by approximately 40%, compared with the computation of the full Jacobian $\boldsymbol{J}5$.

Moreover, numerical experiments presented in section 5 show that the convergence of the Newton iteration is barely affected when using the simplified Jacobian $\boldsymbol{J}3$.

### 3.4. Geometric multigrid method

To efficiently solve the final linear system, assembled from (13) with the simplified Jacobian $\boldsymbol{J}3$, the geometric multigrid method is employed. As described in [15, 21, 22], the multigrid solver can be constructed in a standard framework, provided that the following key components are properly defined: a hierarchy of grids, the restriction and prolongation operators between successive levels, and an appropriate smoothing operator on each grid.

Let $\{\mathcal{G}_l : l = 0, 1, \ldots, L\}$ denote the hierarchy of $L + 1$ grid levels over the computational domain, where $\mathcal{G}_0$ is the finest grid, on which the finite volume WENO discretization and the corresponding linear system are defined. The coarser grids are obtained successively by aggregating neighboring cells from the finer grids. Specifically, each cell $I_{j,l} \in \mathcal{G}_l$ is composed of two adjacent finer-level cells in $\mathcal{G}_{l-1}$ as

$$I_{j,l} = I_{2j-1,l-1} \bigcup I_{2j,l-1}. \tag{15}$$

Moreover, the coarsest grid $\mathcal{G}_L$ usually contains only a few cells, so that the associated linear system can be efficiently solved using a direct method.

On the $l$th grid level, we shall solve a linear system whose $j$th equation is given by

$$\sum_{i=j-1}^{j+1} J_{ji,l}\, \delta U_{i,l} = -R_{j,l}, \tag{16}$$

where $J_{ji,l} = 0$ for $i < 1$ or $i > N_l$, and $N_l$ denotes the number of cells in $\mathcal{G}_l$. At the finest level $l = 0$, the coefficient blocks $J_{ji,0}$ and the right-hand side $R_{j,0}$ are defined as

$$J_{ji,0} = J_{ji} = \alpha \left\| R_j^{(n)} \right\|_{\ell^1} \delta_{ji} \mathbf{I} + \left( \frac{\partial R_j}{\partial \overline{U}_i} \right)^{(n)}, \quad R_{j,0} = R_j^{(n)}.$$

For coarser levels $l \geq 1$, they are constructed using a restriction procedure similar to that in [15, 21], yielding

$$J_{ji,l} = \sum_{\xi=2j-1}^{2j} \sum_{\eta=2i-1}^{2i} J_{\xi\eta,l-1}, \quad R_{j,l} = \sum_{\xi=2j-1}^{2j} \left( R_{\xi,l-1} + \sum_{\eta=\xi-1}^{\xi+1} J_{\xi\eta,l-1} \delta U_{\eta,l-1} \right), \tag{17}$$

where $\delta U_{\eta,l-1}$, for all $\eta$, are the current approximation to the $(l-1)$th-level system. Since the simplified Jacobian $\boldsymbol{J}3$ is employed, the off-diagonal coefficient blocks immediately reduce to $J_{j,j+1,l} = J_{2j,2j+1,l-1}$ and $J_{j,j-1,l} = J_{2j-1,2j-2,l-1}$.

Once the coarse-level approximation $\delta U_{j,l}$ has been obtained, the fine-level approximation is corrected through a simple prolongation operator as

$$\delta U_{i,l-1} \coloneqq \delta U_{i,l-1} + \delta U_{j,l}, \quad \text{for } i = 2j-1, 2j. \tag{18}$$

In addition, to effectively damp out the high-frequency components of the error, a symmetric block SOR method is adopted as the smoother on each grid level. It is implemented as a fast sweeping iteration that performs one forward and one backward sweep over all cells. More precisely, the $j$th component $\delta U_{j,l}$ of the solution to the $l$th-level linear system (16) is updated as

$$\delta U_{j,l} \longleftarrow (1-\omega)\delta U_{j,l} - \omega J_{jj,l}^{-1} \left( R_{j,l} + J_{j,j-1,l}\delta U_{j-1,l} + J_{j,j+1,l}\delta U_{j+1,l} \right), \tag{19}$$

where $\omega$ is the relaxation parameter, and the index $j$ is traversed forward from 1 to $N_l$ and then backward, thus completing one symmetric SOR iteration.

With the components described above, an $(L + 1)$-level multigrid method can be recursively constructed to solve the finest-level linear system. Specifically, this multigrid iteration at level $l$ following a $V$-cycle procedure, denoted as $\delta \boldsymbol{U}_l^{m+1} = \mathrm{MG}_L(l, \delta \boldsymbol{U}_l^m)$, is summarized in **Algorithm 2**, where $\delta \boldsymbol{U}_l^m$ represents the $m$th approximation to the $l$th-level linear system, with its $j$th component given by $\delta U_{j,l}^m$. Accordingly, the multigrid iteration for the finest-level system can be written as $\delta \boldsymbol{U}_0^{m+1} = \mathrm{MG}_L(0, \delta \boldsymbol{U}_0^m)$.

---

**Algorithm 2** Multigrid $V$-cycle iteration for the $l$th-level system: $\delta \boldsymbol{U}_l^{m+1} = \mathrm{MG}_L(l, \delta \boldsymbol{U}_l^m)$

1: **if** $l = L$ **then** /* *Coarsest level* */
2: Solve the corresponding linear system directly to obtain $\delta \boldsymbol{U}_L^{m+1}$;
3: **else**
4: *Pre-smoothing*: Update the $l$th-level approximation from $\delta \boldsymbol{U}_l^m$ to $\delta \widehat{\boldsymbol{U}}_l$ by performing $\nu_1$ steps of the smoother;
5: *Coarse-level correction*:
6: Construct the $(l+1)$th-level coefficient blocks $J_{ji,l+1}$ and right-hand side $R_{j,l+1}$;
7: Initialize the coarse-level approximation as $\delta \boldsymbol{U}_{l+1}^0 = 0$;
8: Recursively call the multigrid iteration at level $l+1$ to obtain $\delta \widetilde{\boldsymbol{U}}_{l+1} = \mathrm{MG}_L(l+1, \delta \boldsymbol{U}_{l+1}^0)$;
9: Update the $l$th-level approximation by prolongation

$$\delta \widetilde{U}_{i,l} \longleftarrow \delta \widehat{U}_{i,l} + \delta \widetilde{U}_{j,l+1} \quad \text{for all } j \text{ and } i = 2j-1, 2j;$$

10: *Post-smoothing*: Update the $l$th-level approximation from $\delta \widetilde{\boldsymbol{U}}_l$ to $\delta \boldsymbol{U}_l^{m+1}$ by performing $\nu_2$ steps of the smoother;
11: **end if**

---

By incorporating the multigrid iteration for the linear system into the framework of the Newton-based solver described in **Algorithm 1**, a third-order NMGM can be naturally developed for solving the underlying nonlinear problem arising from the finite volume WENO discretization on the target grid. However, it is worth noting that we have not yet specified how to determine an appropriate initial guess for the iteration, which plays a crucial role in ensuring the robustness and convergence of Newton's method.

To address this issue, we adopt a hierarchical initialization strategy, which is conceptually similar to the philosophy of the full multigrid approach [22, 23]. Specifically, we first solve the nonlinear problem resulting from the discretization on the coarsest grid. At this grid level, because of the small number of unknowns, Newton's method, with the linearized system solved directly, can be efficiently applied and converge readily with a moderate choice of the initial guess. On finer grid levels, the corresponding nonlinear problem is solved using the NMGM, with the initial guess obtained by cell averaging a high-order polynomial interpolated from the coarser-level solution. This procedure allows for efficient and accurate construction of initial approximations on finer grids, and ultimately yields the complete NMGM solver, as summarized in **Algorithm 3**. Additionally, the computation of the initial cell averages from the coarser-level solution is detailed in **Appendix A**.

**Remark 1.** It is apparent that the entire procedure described above can also be applied to the unsimplified regularized linear system, in which the full Jacobian $\boldsymbol{J}5$, rather than the simplified $\boldsymbol{J}3$, is retained in the coefficient matrix. To distinguish between the two formulations throughout the paper, the NMGM solvers based on the simplified and full Jacobians are referred to as NMGM-3 and NMGM-5, respectively.

## 4. Extension of the NMGM to two-dimensional problems

The third-order NMGM introduced in the previous section can be readily extended to the two-dimensional problem (2) on rectangular grids. Since the overall algorithmic framework remains essentially unchanged, this section focuses on the key modifications and distinctive features that arise in the two-dimensional formulation.

### *4.1. Dimension-by-dimension WENO reconstruction and Jacobian simplification*

Let the computational domain be uniformly partitioned into a rectangular grid of $N_x \times N_y$ cells, where the $(i, j)$th cell is defined as $I_{ij} = \left[x_{i-\frac{1}{2}}, x_{i+\frac{1}{2}}\right] \times \left[y_{j-\frac{1}{2}}, y_{j+\frac{1}{2}}\right]$, $i = 1, 2, \ldots, N_x$, $j = 1, 2, \ldots, N_y$, and the cell edge lengths are $\Delta x = x_{i+\frac{1}{2}} - x_{i-\frac{1}{2}}$ and $\Delta y = y_{j+\frac{1}{2}} - y_{j-\frac{1}{2}}$ in the $x$- and $y$-directions, respectively. The finite volume formulation of (2) over $I_{ij}$ yields

$$R_{ij}(\overline{\boldsymbol{U}}_{ij}) \triangleq \frac{\left(\hat{f}_{i+\frac{1}{2},j} - \hat{f}_{i-\frac{1}{2},j}\right)}{\Delta x} + \frac{\left(\hat{g}_{i,j+\frac{1}{2}} - \hat{g}_{i,j-\frac{1}{2}}\right)}{\Delta y} - \frac{S_{ij}}{\Delta x \Delta y} = 0, \tag{20}$$

**Algorithm 3** Newton multigrid (NMGM) solver

1: Given the tolerance *Tol* and the maximum number of Newton iterations $n_{\max}$;
2: Initialize the coarsest-level cell averages using the given initial guess $\overline{\boldsymbol{U}}_L^{(0)}$, and set $n = 0$;
3: **repeat** /* *Newton iteration on the coarsest grid level* */
4: Solve the coarsest-level regularized linear system (13) directly to obtain the increment $\delta\boldsymbol{U}_L^{(n)}$;
5: Update the cell averages according to (14), which in vector form reads $\overline{\boldsymbol{U}}_L^{(n+1)} = \overline{\boldsymbol{U}}_L^{(n)} + \tau\,\delta\boldsymbol{U}_L^{(n)}$;
6: Evaluate the global residual at this level by $res = \frac{1}{N_L}\sum_{j=1}^{N_L}\left\|R_{j,L}^{(n+1)}\right\|_{\ell^1}$, and set $n \leftarrow n+1$;
7: **until** $res < Tol$ or $n \geq n_{\max}$.
8: **for** $l$ from $L-1$ **to** 0 **do** /* *Solver on the lth grid level* */
9: Evaluate the initial cell-averaged approximation $\overline{\boldsymbol{U}}_l^{(0)}$ by high-order interpolation from the coarser-level solution $\overline{\boldsymbol{U}}_{l+1}^{(n)}$;
10: Set $n = 0$;
11: **repeat** /* *Newton multigrid iteration on the lth grid level* */
12: Initialize $\delta\boldsymbol{U}_l^0 = 0$;
13: **for** $m$ from 0 to $\nu_3 - 1$ **do** /* *Solve the lth-level regularized system* (13) *by* $\nu_3$ *multigrid iterations* */
14: Perform one $(L-l+1)$-level multigrid iteration to obtain $\delta\boldsymbol{U}_l^{m+1} = \mathrm{MG}_{L-l+1}(0, \delta\boldsymbol{U}_l^m)$;
15: **end for**
16: Update the cell averages according to (14), which in vector form reads $\overline{\boldsymbol{U}}_l^{(n+1)} = \overline{\boldsymbol{U}}_l^{(n)} + \tau\,\delta\boldsymbol{U}_l^{\nu_3}$;
17: Evaluate the global residual at this level by $res = \frac{1}{N_l}\sum_{j=1}^{N_l}\left\|R_{j,l}^{(n+1)}\right\|_{\ell^1}$, and set $n \leftarrow n+1$;
18: **until** $res < Tol$ or $n \geq n_{\max}$.
19: **end for**

where $\overline{\boldsymbol{U}}_{ij} = \left\{\overline{U}_{pq}\right\}_{(p,q)\in\mathcal{T}_{ij}}$ is the collection of all cell averages within the stencil $\mathcal{T}_{ij}$ used to evaluate the local residual $R_{ij}$. The symbols $\hat{f}_{i\pm\frac12,j}$ and $\hat{g}_{i,j\pm\frac12}$ represent the edge-averaged numerical fluxes across the vertical and horizontal interfaces of the cell, respectively, and $S_{ij}$ denotes the cell-integrated source term over $I_{ij}$. Following [16], we employ a dimension-by-dimension approach for the evaluation of these quantities.

To achieve third-order accuracy, the numerical fluxes $\hat{f}_{i\pm\frac12,j}$ and $\hat{g}_{i,j\pm\frac12}$ are computed using the two-point Gauss-Legendre quadrature rule, as in (9). For instance, we have

$$\hat{f}_{i+\frac12,j} = \frac{1}{\Delta y}\int_{y_{j-\frac12}}^{y_{j+\frac12}} F\left(U(x_{i+\frac12}, y)\right)\mathrm{d}y \approx \frac12\sum_{\beta=1}^{2}\hat{F}\left(U_{x^-_{i+\frac12},y_j^\beta}, U_{x^+_{i+\frac12},y_j^\beta}\right),$$

where $y_j^\beta$, $\beta = 1,2$, are the Gauss-Legendre quadrature points in $\left[y_{j-\frac12}, y_{j+\frac12}\right]$, and $U_{x^\pm_{i+\frac12},y_j^\beta} = U^\pm_{h,i+\frac12}(y_j^\beta)$ are the pointwise approximations to $U(x_{i+\frac12}, y_j^\beta)$. Here $U^\pm_{h,i+\frac12}(y)$ are polynomials obtained by interpolating the interface values $U^\pm_{i+\frac12,j-1}$, $U^\pm_{i+\frac12,j}$, and $U^\pm_{i+\frac12,j+1}$, which are generated from the one-dimensional WENO reconstruction in the $x$-direction.

For the cell-integrated source term $S_{ij}$, it is straightforward to see that its first component is zero. For the second component $S_{2,ij}$, we compute it by successively applying the two-point Gauss-Legendre quadrature rule in the $y$- and $x$-directions as follows

$$S_{2,ij} = \int_{I_{ij}} -ghb_x\,\mathrm{d}x\,\mathrm{d}y \approx -\frac{\Delta y}{2}\sum_{\beta=1}^{2}\int_{x_{i-\frac12}}^{x_{i+\frac12}} gh_h^\beta(x)b_x(x, y_j^\beta)\,\mathrm{d}x \approx -\frac{\Delta x\Delta y}{4}\sum_{\beta=1}^{2}\sum_{\gamma=1}^{2} gh_h^\beta(x_i^\gamma)b_x(x_i^\gamma, y_j^\beta),$$

where $x_i^\gamma$, $\gamma = 1,2$, are the Gauss-Legendre quadrature points in $\left[x_{i-\frac12}, x_{i+\frac12}\right]$, and $h_h^\beta(x)$ are the interpolating polynomials constructed from the interface values along $y = y_j^\beta$, as in the one-dimensional case. Similar to [16], the third component $S_{3,ij}$ is approximated by

$$S_{3,ij} = \int_{I_{ij}} -ghb_y\,\mathrm{d}x\,\mathrm{d}y \approx -\frac{\Delta x}{2}\sum_{\gamma=1}^{2}\int_{y_{j-\frac12}}^{y_{j+\frac12}} g\tilde{h}_h^\gamma(y)b_y(x_i^\gamma, y)\,\mathrm{d}y \approx -\frac{\Delta x\Delta y}{4}\sum_{\gamma=1}^{2}\sum_{\beta=1}^{2} g\tilde{h}_h^\gamma(y_j^\beta)b_y(x_i^\gamma, y_j^\beta),$$

where $\tilde{h}_h^\gamma(y)$ denote the interpolating polynomials constructed along $x = x_i^\gamma$.

From the above formulation, we can identify that the stencil $\mathcal{T}_{ij}$ required to evaluate the local residual $R_{ij}$ consists of 21 cells, as illustrated in Figure 3a. Consequently, the full Jacobian matrix arising from the regularized Newton iteration applied to (20) has a banded block structure, with at most 21 nonzero blocks per row, each of size $3 \times 3$. Similar to the one-dimensional case, we denote this matrix by $\boldsymbol{J}21$ and examine the entrywise $\ell^1$-norms of all its nonzero blocks. Typical distributions of these block magnitudes in representative rows are displayed in Figures 4a and 4b, from which we observe clearly distinguishable magnitudes among the blocks.

For convenience, let us introduce a reduced stencil $\hat{\mathcal{T}}_{ij}$, consisting of 9 cells, whose configuration is shown in Figure 3b. Then we find that the blocks associated with the cells in $\hat{\mathcal{T}}_{ij}$ are generally of considerably larger magnitude than those associated with the cells in $\mathcal{T}_{ij} \setminus \hat{\mathcal{T}}_{ij}$, i.e., the cells in $\mathcal{T}_{ij}$ but not in $\hat{\mathcal{T}}_{ij}$. More specifically, for the full $\boldsymbol{J}21$ obtained from the numerical experiment in subsection 5.2.1 on a $16 \times 8$ grid, the magnitudes of blocks corresponding to the cells in $\hat{\mathcal{T}}_{ij}$ consistently exceed 0.5, with nearly 79% greater than 1, whereas all remaining blocks have magnitudes below 0.5, with about 83% even smaller than $1 \times 10^{-6}$. For the numerical experiment in subsection 5.2.2, these two fractions are 100% and 98%, respectively.

Motivated by these observations, we approximate $\boldsymbol{J}21$ by retaining only the blocks associated with the cells in $\hat{\mathcal{T}}_{ij}$ and discarding the rest. The resulting simplified Jacobian is denoted by $\boldsymbol{J}9$, whose pattern of nonzero blocks is shown in Figure 4c. With such a simplification, the computational cost of evaluating the Jacobian is reduced by nearly 57%.

Denote by $J^{ij}_{pq}$ the $(p, q)$th Jacobian block associated with the $(i, j)$th residual equation (20) in the regularized Newton iteration, that is, $J^{ij}_{pq} = \alpha \left\| R^{(n)}_{ij} \right\|_{\ell^1} \delta_{ip}\delta_{jq}\mathbf{I} + \left( \partial R_{ij} / \partial \overline{U}_{pq} \right)^{(n)}$. Then the modified linearization of (20), corresponding to the simplified $\boldsymbol{J}9$, can be written as

$$\sum_{(p,q)\in\hat{\mathcal{T}}_{ij}} J^{ij}_{pq}\,\delta U^{(n)}_{pq} = J^{ij}_{ij}\,\delta U^{(n)}_{ij} + \sum_{p=i-2,\,p\neq i}^{i+2} J^{ij}_{pj}\,\delta U^{(n)}_{pj} + \sum_{q=j-2,\,q\neq j}^{j+2} J^{ij}_{iq}\,\delta U^{(n)}_{iq} = -R^{(n)}_{ij}. \tag{21}$$

Additionally, the linearization corresponding to $\boldsymbol{J}21$ can be obtained from the above equation by replacing the index set $\hat{\mathcal{T}}_{ij}$ with the full stencil $\mathcal{T}_{ij}$.

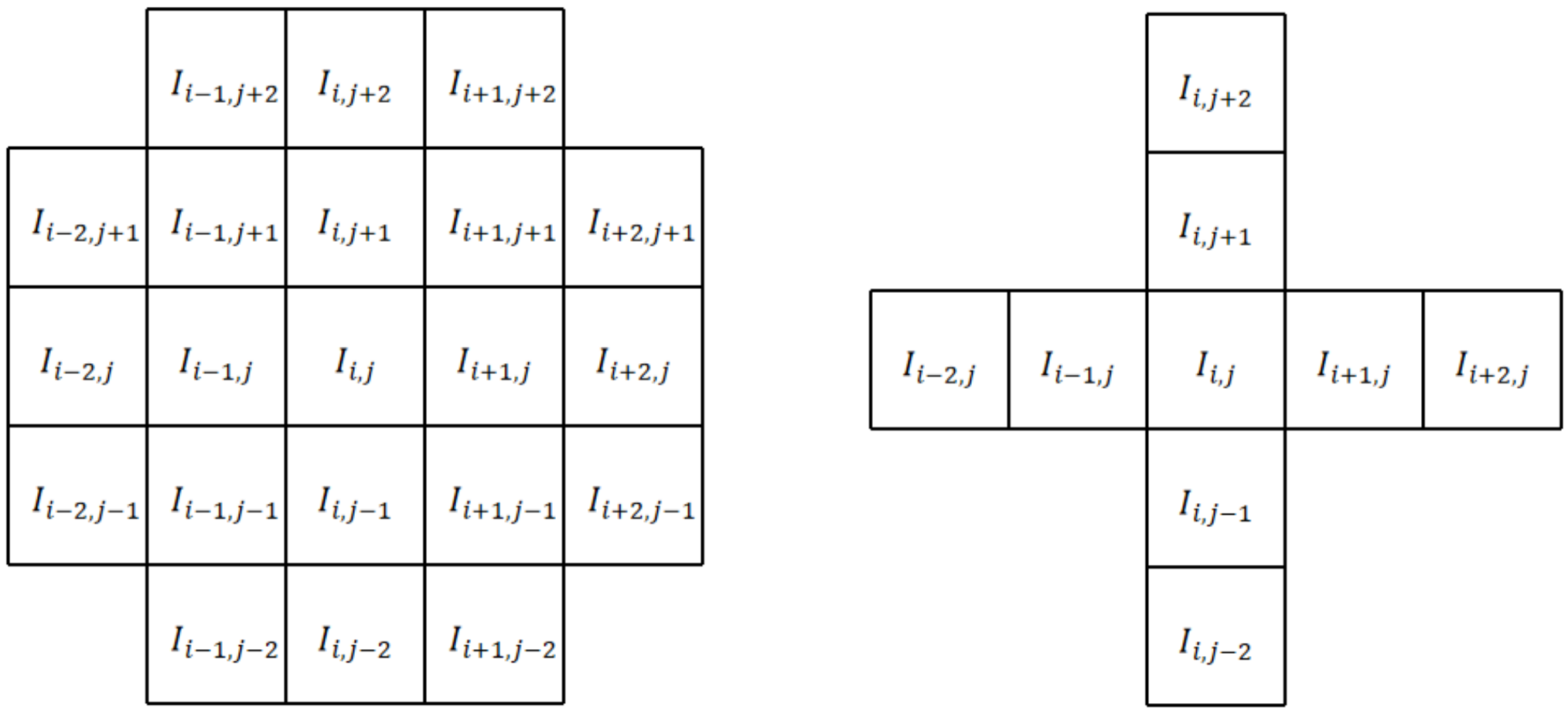

(a) Stencil $\mathcal{T}_{ij}$ for the discretization.

(b) Nonzero block stencil $\hat{\mathcal{T}}_{ij}$ of the simplified Jacobian.

Figure 3: Stencil configurations for the third-order discretization and the nonzero blocks of the simplified Jacobian associated with $I_{ij}$.

## 4.2. *Multigrid smoothing and restriction operators*

As a natural extension of the one-dimensional smoother, we employ the block SOR method combined with a prescribed sweeping strategy as the smoother in the multigrid algorithm on rectangular grids. For the linear system (21), it follows that the

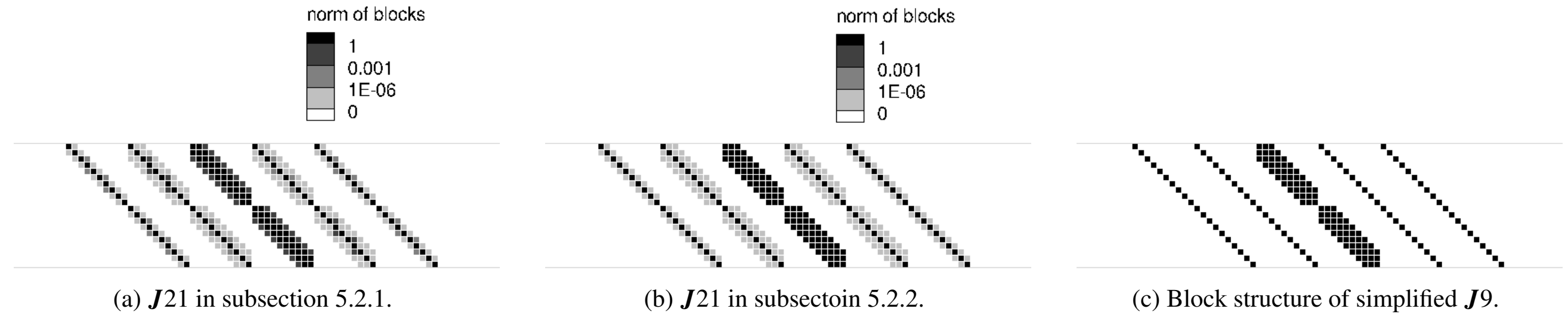

(a) $\boldsymbol{J}21$ in subsection 5.2.1. (b) $\boldsymbol{J}21$ in subsectoin 5.2.2. (c) Block structure of simplified $\boldsymbol{J}9$.

Figure 4: Representative rows of the Jacobian matrices for a $10 \times 10$ grid: typical block magnitudes for $\boldsymbol{J}21$, and block structure of the simplified $\boldsymbol{J}9$.

update of the $(i, j)$th component $\delta U_{ij}^{(n)}$ takes the form

$$\delta U_{ij}^{(n)} \longleftarrow (1-\omega)\,\delta U_{ij}^{(n)} - \omega \left(J_{ij}^{ij}\right)^{-1} \left( R_{ij}^{(n)} + \sum_{p=i-2,\, p\neq i}^{i+2} J_{pj}^{ij}\,\delta U_{pj}^{(n)} + \sum_{q=j-2,\, q\neq j}^{j+2} J_{iq}^{ij}\,\delta U_{iq}^{(n)} \right). \tag{22}$$

This update is carried out according to a fast sweeping strategy, where each iteration consists of four-direction alternating sweeps described by

$$\begin{aligned}
&(\mathrm{D1}): \quad i = 1, 2, \ldots, N_x \text{ (outer loop)}; \qquad j = 1, 2, \ldots, N_y \text{ (inner loop)};\\
&(\mathrm{D2}): \quad i = N_x, \ldots, 2, 1 \text{ (outer loop)}; \qquad j = 1, 2, \ldots, N_y \text{ (inner loop)};\\
&(\mathrm{D3}): \quad i = N_x, \ldots, 2, 1 \text{ (outer loop)}; \qquad j = N_y, \ldots, 2, 1 \text{ (inner loop)};\\
&(\mathrm{D4}): \quad i = 1, 2, \ldots, N_x \text{ (outer loop)}; \qquad j = N_y, \ldots, 2, 1 \text{ (inner loop)}.
\end{aligned}$$

Indeed, such a sweeping procedure is a direct generalization of the one-dimensional symmetric sweeping strategy. Moreover, as explained in [24, 25], it captures the dominant propagation directions inherited from the hyperbolic structure of the governing equations, thereby enhancing both smoothing efficiency and robustness.

To construct the restriction and prolongation operators, a hierarchy of grids must first be specified. As in the one-dimensional case, this hierarchy of grids is generated by successively aggregating neighboring cells on the finer grid level. It follows that the $(i, j)$th coarse-level cell $I_{ij,l} \in \mathcal{G}_l$ is composed of four adjacent finer-level cells in $\mathcal{G}_{l-1}$, that is,

$$I_{ij,l} = I_{2i-1,2j-1,l-1} \bigcup I_{2i,2j-1,l-1} \bigcup I_{2i-1,2j,l-1} \bigcup I_{2i,2j,l-1}. \tag{23}$$

With this geometric relationship, both restriction and prolongation operators can be locally implemented between adjacent grid levels.

The restriction operator is applied to obtain the coarse-level linear system by assembling the coefficient blocks $J_{pq,l}^{ij}$ and the right-hand side $R_{ij,l}$ from the corresponding quantities on the $(l-1)$th grid level. Analogous to the one-dimensional construction (17), its two-dimensional counterpart on rectangular grids gives rise to

$$J_{pq,l}^{ij} = \sum_{\imath=2i-1}^{2i} \sum_{\jmath=2j-1}^{2j} \sum_{\xi=2p-1}^{2p} \sum_{\eta=2q-1}^{2q} J_{\xi\eta,l-1}^{\imath\jmath}, \qquad R_{ij,l} = \sum_{\imath=2i-1}^{2i} \sum_{\jmath=2j-1}^{2j} \left( R_{\imath\jmath,l-1} + \sum_{(\xi,\eta)\in\hat{\mathcal{T}}_{\imath\jmath,l-1}} J_{\xi\eta,l-1}^{\imath\jmath}\,\delta U_{\xi\eta,l-1} \right), \tag{24}$$

where $\hat{\mathcal{T}}_{\imath\jmath,l-1}$ denotes the stencil of nonzero coefficient blocks associated with the $(\imath, \jmath)$th equation on the $(l-1)$th grid level. At the finest level $l = 0$, we set $J^{ij}_{pq,0} = J^{ij}_{pq}$, $R_{ij,0} = R^{(n)}_{ij}$, and $\hat{\mathcal{T}}_{ij,0} = \hat{\mathcal{T}}_{ij}$. It follows that

$$J^{ij}_{ij,1} = \sum_{\xi=2i-1}^{2i} \sum_{\eta=2j-1}^{2j} \left( J^{\xi\eta}_{\xi\eta,0} + J^{\xi\eta}_{\xi,4j-1-\eta,0} + J^{\xi\eta}_{4i-1-\xi,\eta,0} \right), \quad J^{ij}_{i+1,j,1} = \sum_{\eta=2j-1}^{2j} \left( J^{2i-1,\eta}_{2i+1,\eta,0} + \sum_{\xi=2i+1}^{2i+2} J^{2i,\eta}_{\xi\eta,0} \right),$$

$$J^{ij}_{i-1,j,1} = \sum_{\eta=2j-1}^{2j} \left( J^{2i,\eta}_{2i-2,\eta,0} + \sum_{\xi=2i-3}^{2i-2} J^{2i-1,\eta}_{\xi\eta,0} \right), \quad J^{ij}_{i,j+1,1} = \sum_{\xi=2i-1}^{2i} \left( J^{\xi,2j-1}_{\xi,2j+1,0} + \sum_{\eta=2j+1}^{2j+2} J^{\xi,2j}_{\xi\eta,0} \right), \quad J^{ij}_{i,j-1,1} = \sum_{\xi=2i-1}^{2i} \left( J^{\xi,2j}_{\xi,2j-2,0} + \sum_{\eta=2j-3}^{2j-2} J^{\xi,2j-1}_{\xi\eta,0} \right),$$

and $J^{ij}_{pq,1} = 0$ for all other index pairs $(p, q)$. For coarser levels $l \geq 2$, we further obtain

$$J^{ij}_{ij,l} = \sum_{\xi=2i-1}^{2i} \sum_{\eta=2j-1}^{2j} \left( J^{\xi\eta}_{\xi\eta,l-1} + J^{\xi\eta}_{\xi,4j-1-\eta,l-1} + J^{\xi\eta}_{4i-1-\xi,\eta,l-1} \right), \quad J^{ij}_{i+1,j,l} = J^{2i,2j-1}_{2i+1,2j-1,l-1} + J^{2i,2j}_{2i+1,2j,l-1},$$

$$J^{ij}_{i-1,j,l} = J^{2i-1,2j-1}_{2i-2,2j-1,l-1} + J^{2i-1,2j}_{2i-2,2j,l-1}, \quad J^{ij}_{i,j+1,l} = J^{2i-1,2j}_{2i-1,2j+1,l-1} + J^{2i,2j}_{2i,2j+1,l-1}, \quad J^{ij}_{i,j-1,l} = J^{2i-1,2j-1}_{2i-1,2j-2,l-1} + J^{2i,2j-1}_{2i,2j-2,l-1},$$

with all other blocks $J^{ij}_{pq,l}$ vanishing. This implies that, for all levels $l \geq 1$, the stencil $\hat{\mathcal{T}}_{ij,l}$ involves at most five cells, thus forming a block-pentadiagonal linear system at each coarse grid level. Consequently, the present multigrid algorithm can be efficiently implemented with much lower per-iteration computational cost than its counterpart based on the full Jacobian $\boldsymbol{J}21$, which results in denser coarse-level stencils.

As for the prolongation operator, the two-dimensional analogue of (18) is given by

$$\delta U_{ij,l-1} \coloneqq \delta U_{ij,l-1} + \delta U_{\xi\eta,l}, \quad \text{for} \;\; i = 2\xi - 1, 2\xi \;\; \text{and} \;\; j = 2\eta - 1, 2\eta. \tag{25}$$

Now, the two-dimensional NMGM solver can be formulated using the above components. Following the one-dimensional notation, the solver based on the simplified $\boldsymbol{J}9$ is referred to as NMGM-9, while the version using the full $\boldsymbol{J}21$ is referred to as NMGM-21.

## 5. Numerical experiments

Several numerical experiments are conducted in this section to demonstrate the accuracy, efficiency, and robustness of the proposed third-order NMGM under a nondimensional framework for both one- and two-dimensional SWEs. Unless otherwise specified, the following solver settings are adopted. The parameter $\alpha$ used to formulate the regularized linear system (13) is set to 3, and the damping parameter $\tau$ in the Newton update (14) is chosen as 0.6. For the one-dimensional experiments, the number of multigrid iterations per Newton step is taken as $\nu_3 = 2$. Within each multigrid iteration, the SOR relaxation parameter is selected as $\omega = 1.0$, the numbers of pre- and post-smoothing steps are set to $\nu_1 = \nu_2 = 2$, and the total number of grid levels is fixed at $L + 1 = 4$. Following [15, 4], the HLL flux is employed, since it consistently outperforms the LLF flux. However, as also reported in [15], the HLL flux may fail to converge in two-dimensional settings. Therefore, we adopt the LLF flux in the two-dimensional experiments, with appropriately selected parameters. Specifically, we set $\nu_3 = 3$, $\omega = 0.1$, and determine the total number of grid levels such that the coarsest grid consists of $8 \times 4$ cells, while the smoothing steps $\nu_1$ and $\nu_2$ are chosen in a grid-dependent manner to be specified later.

All simulations are implemented in `C++` and executed on the High Performance Computing Platform of Nanjing University of Aeronautics and Astronautics, China. They are terminated once the global residual norm *res* drops below a tolerance of $Tol = 1 \times 10^{-11}$ for both the one- and two-dimensional cases.

### 5.1. One-dimensional case

In the one-dimensional case, numerical tests are primarily performed on uniform grids with $N = 96$, 192, 384, and 768 cells. On these grids, it is observed that the numerical differentiation (12) with a too small perturbation $\epsilon$, such as $\epsilon = 10^{-6}$ used for the first-order discretization in [15], may lead to numerical instabilities and result in NaN values. To ensure stable computations, we empirically adopt $\epsilon = 0.2$, 0.15, 0.1, and 0.05 for these four grid resolutions, respectively. For finer grids, setting $\epsilon = 0.01$ is usually sufficient. It is also observed that, once numerical stability is ensured, these specific choices of $\epsilon$ do not significantly affect the convergence behavior of the present NMGM solver.

*5.1.1. Smooth subcritical flow*

The first experiment considers a smooth subcritical flow in a river of constant width $\sigma = 1$ with a non-flat bottom topography defined as

$$b(x) = 0.2\exp\left(-\frac{(x+1)^2}{2}\right) + 0.3\exp\left(-(x-1.5)^2\right), \quad x \in [-10,\ 10]. \tag{26}$$

The boundary conditions $h = 1$ and $hu = 1$ are applied at $x = \pm 10$. This problem has been studied in [15, 26], where the exact solution is shown to satisfy

$$u^3 + (2gb - 2g - 1)\,u + 2g = 0, \quad hu = 1. \tag{27}$$

We solve the above algebraic system to obtain the exact solution for validation and comparison.

The $\ell^1$ errors for the cell averages of the conservative variables $h$ and $hu$, together with the corresponding numerical orders of accuracy, are reported in Table 1. The results show a steady error reduction under grid refinement, with convergence rates approaching third order, which is consistent with the expected third-order accuracy of the present finite volume WENO discretization.

Table 1: Smooth subcritical flow: $\ell^1$ errors and numerical orders of accuracy.

| No. of cells | $h$ | | $hu$ | |
|---|---|---|---|---|
| | $\ell^1$ error | order | $\ell^1$ error | order |
| 96 | 4.64E-4 | - | 1.05E-3 | - |
| 192 | 6.58E-5 | 2.82 | 1.60E-4 | 2.71 |
| 384 | 9.23E-6 | 2.83 | 2.12E-5 | 2.92 |
| 768 | 1.18E-6 | 2.97 | 2.42E-6 | 3.13 |

Moreover, it is observed from Table 1 that the third-order NMGM solver attains $\ell^1$ errors on the order of $10^{-3}$-$10^{-4}$ using a grid of 96 cells. By comparison, the first-order NMGM solver introduced in [15] requires substantially finer grids to achieve comparable accuracy. This difference is illustrated in Figure 5, which compares the steady-state profiles of the surface level $h + b$ and the velocity $u$ obtained on the 96-cell grid with those computed by the first-order method on grids of 768 and 1536 cells, along with the exact solutions. It can be seen that, even on the 1536-cell grid, the first-order solution still exhibits visible discrepancies, whereas accurate profiles are already obtained on the much coarser 96-cell grid with the third-order discretization. Further numerical tests on finer grids indicate that the first-order method requires grids with several thousand cells to reach accuracy levels comparable to the third-order solution on the 96-cell grid, at a significantly higher computational cost. These results highlight the advantage of the proposed approach in achieving the desired accuracy on much coarser grids.

To investigate the convergence behavior and computational cost of the proposed NMGM solver, we list the number of Newton iterations and the CPU times for NMGM-3 and NMGM-5 on a sequence of refined grids in Table 2. For reference, the corresponding results obtained with a Newton-SOR solver are also included. Here, the Newton-SOR solver refers to a single-level variant of the NMGM framework, in which the inner multigrid iteration is replaced by a symmetric block SOR method. Within each Newton iteration, this inner SOR iteration is terminated when either the $\ell^1$-norm of the variation between two successive iterates falls below $1 \times 10^{-7}$ or the maximum number of 10 iterations is reached. In addition, for computations on each grid, the initial guess is obtained by interpolating the steady-state solution from the preceding coarser grid, as described in **Appendix A**. As for the reported CPU times, they measure only the computational cost on the associated grid level.

The results show that, by incorporating the multigrid strategy, the NMGM solvers significantly reduce both the number of Newton iterations and the CPU time when compared with the reference Newton-SOR solver on the tested grids. Within the NMGM framework, NMGM-3 requires exactly the same number of Newton iterations as NMGM-5, indicating that the use of the simplified Jacobian $\boldsymbol{J}3$ does not affect the convergence of the Newton iteration. At the same time, NMGM-3 consistently achieves lower CPU times than NMGM-5, with a reduction of approximately 30-35% across all tests. As discussed in subsection 3.3, this improvement can be attributed to the reduced computational cost per Newton iteration resulting from the simplified Jacobian, while preserving the robustness of the NMGM strategy.

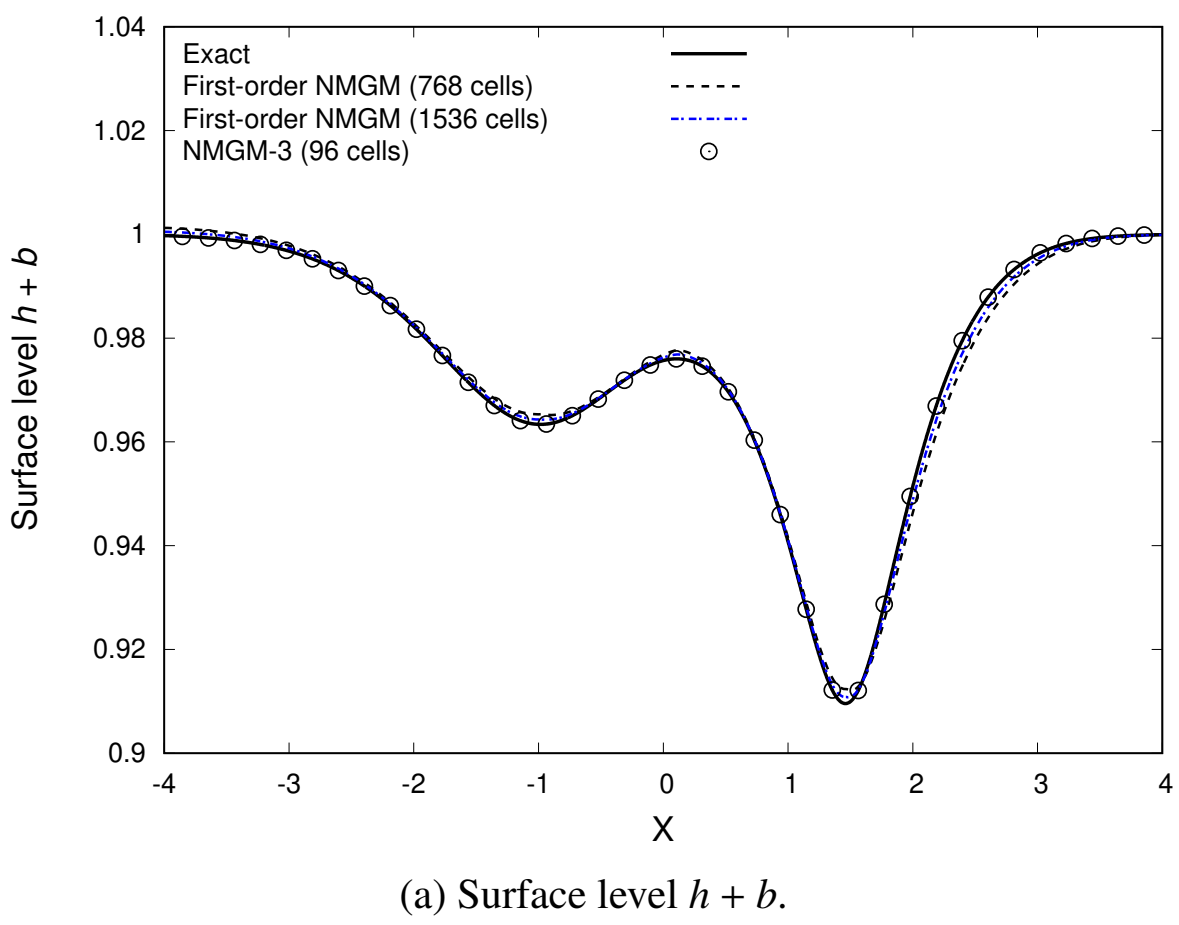

(a) Surface level $h + b$.

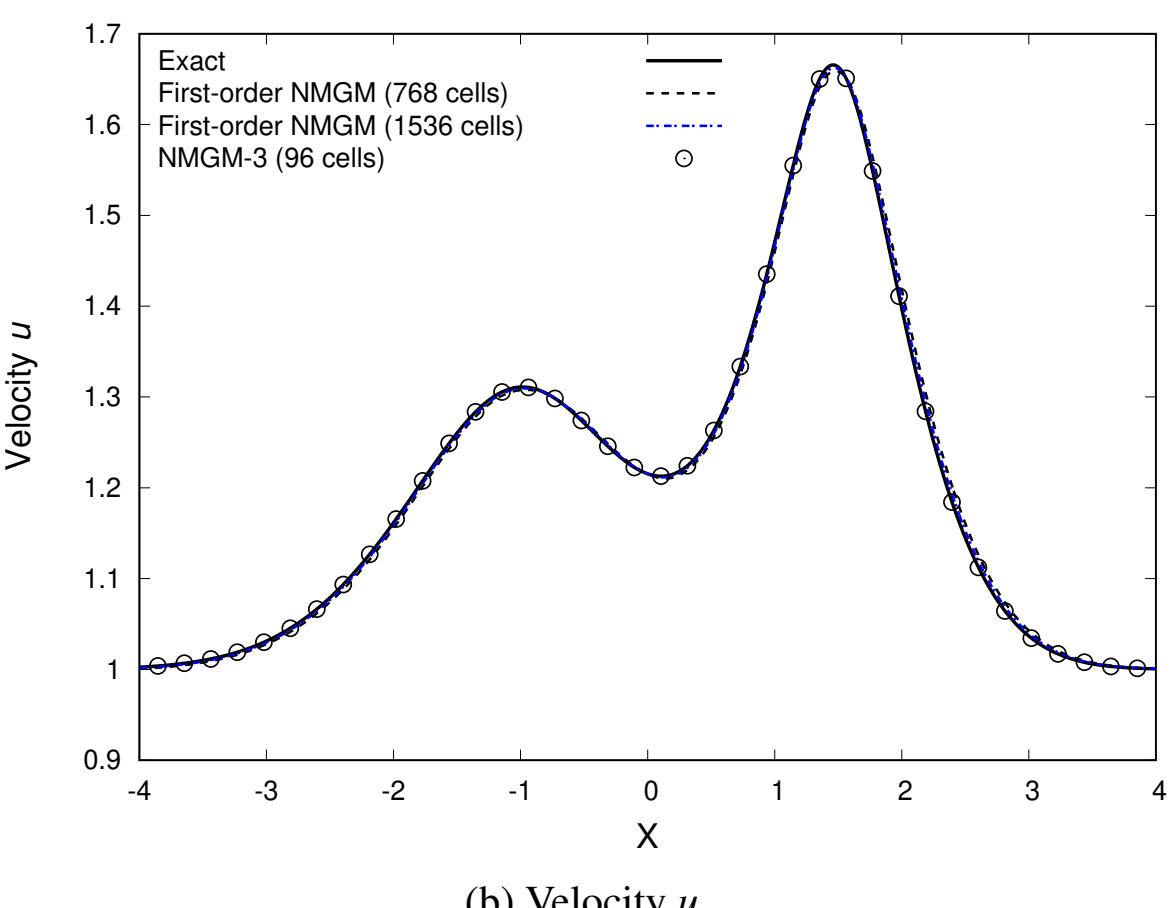

(b) Velocity $u$.

Figure 5: Smooth subcritical flow: comparison of steady-state solutions obtained with the first- and third-order discretizations.

Table 2: Smooth subcritical flow: efficiency comparison of the proposed NMGM solvers and a reference Newton-SOR method.

| No. of cells | No. of Newton iterations | | | CPU time (s) | | |
|---|---|---|---|---|---|---|
| | Newton-SOR | NMGM-5 | NMGM-3 | Newton-SOR | NMGM-5 | NMGM-3 |
| 96 | 116 | 67 | 67 | 22 | 17 | 11 |
| 192 | 110 | 66 | 66 | 44 | 36 | 23 |
| 384 | 87 | 58 | 58 | 80 | 70 | 46 |
| 768 | 79 | 48 | 48 | 182 | 140 | 93 |

#### *5.1.2. Subcritical and transcritical flows over a bump with varying channel width*

The second experiment investigates steady-state flows over a bump in a channel of length 25 with varying width. This problem is also a classical benchmark that has been widely used in the literature to assess the accuracy and robustness of numerical schemes for different flow regimes and channel configurations, see, e.g., [7, 17, 27]. Following the setup in [7], the computational domain is $[0, 25]$, and the bottom topography is given by

$$b(x) = \begin{cases} 0.2 - 0.05(x-10)^2, & x \in [8, 12], \\ 0, & \text{otherwise.} \end{cases} \tag{28}$$

The channel width is either constant, namely $\sigma(x) \equiv 1$, or features a smooth contraction described by

$$\sigma(x) = 1 - \sigma_0 \left( 1 + \cos\left( 2\pi \frac{x - (x_l + x_r)/2}{x_r - x_l} \right) \right), \qquad x \in [x_l, x_r],$$

with $\sigma(x) = 1$ elsewhere. Two contraction configurations are considered: a left-shifted contraction with $[x_l, x_r] = [3.75, 16.25]$, and a right-shifted contraction with $[x_l, x_r] = [8.75, 21.25]$. Depending on the imposed boundary conditions, the flow can be subcritical or transcritical. For the subcritical case, a discharge condition $hu = 4.42$ is prescribed at the upstream boundary $x = 0$, while a water depth condition $h = 2$ is imposed at the downstream boundary $x = 25$. For the transcritical case, the upstream discharge is set to $hu = 1.53$ and the downstream water depth to $h = 0.66$. In the contraction region, the parameter $\sigma_0$ is taken as 0.05 for the subcritical flow and 0.15 for the transcritical flow.

As described in the previous example, the initial guess on each grid level is obtained by interpolating the steady-state solution from the preceding coarser grid. On the coarsest tested grid consisting of 48 cells, the initial guess is constructed according to the state

$$h(x) = 0.5 - b(x), \quad u(x) = 0.$$

*(I) Subcritical flow*

Figure 6 shows the steady-state surface level $h + b$ for the subcritical flow computed by NMGM-3 on a 96-cell grid under three channel width configurations, together with the reference analytical solution from [27] and the bottom topography. For all channel configurations, the numerical results agree well with the reference solution, demonstrating that NMGM-3 can reliably compute steady-state subcritical flows in such settings.

Table 3 further reports the number of Newton iterations required to reach the steady state on four successively refined grids. Since NMGM-3 and NMGM-5 again exhibit nearly identical iteration numbers in these tests, only the results of NMGM-3 are presented. For all configurations, the solver converges robustly with a well-controlled number of Newton iterations. For channels with left- or right-shifted contractions, the iteration numbers generally decrease as the grid is refined, indicating improved efficiency of NMGM-3 on finer grids. By contrast, the constant-width configuration requires more iterations, and the solver shows weak dependence on the grid size. This behavior may be related to the different flow structure in the constant-width case, which produces a flatter steady-state profile and may influence the convergence of the Newton iteration.

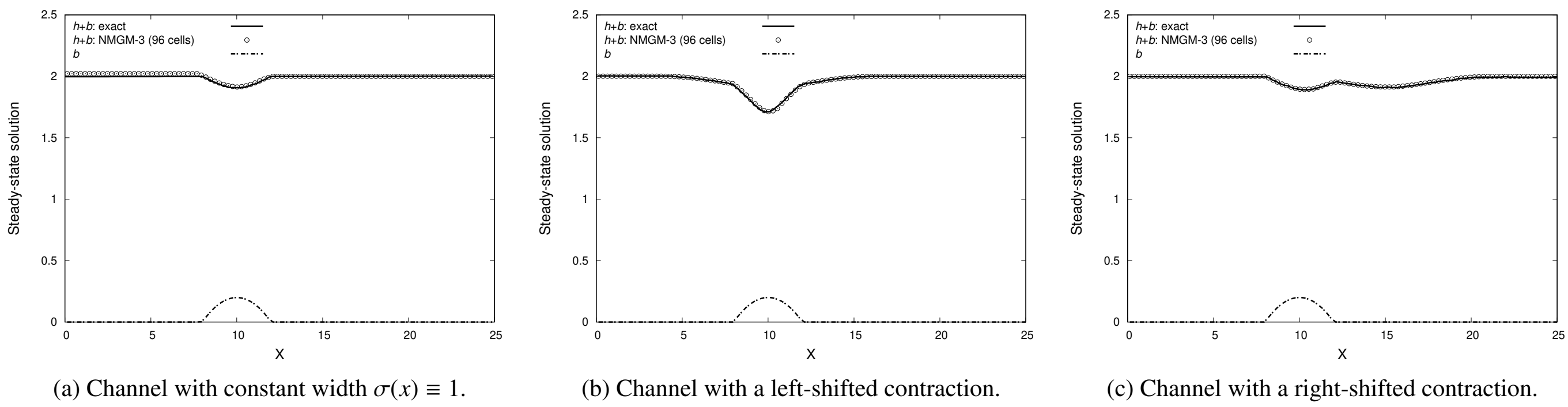

(a) Channel with constant width $\sigma(x) \equiv 1$. (b) Channel with a left-shifted contraction. (c) Channel with a right-shifted contraction.

Figure 6: Subcritical flow over a bump: steady-state surface level $h + b$ obtained by NMGM-3 on a 96-cell grid, compared with the exact solution from [27].

Table 3: Subcritical flow over a bump: number of Newton iterations required by NMGM-3 for different channel widths on successively refined grids.

| No. of cells | Constant width | Left-shifted contraction | Right-shifted contraction |
|---|---|---|---|
| 96 | 152 | 153 | 86 |
| 192 | 129 | 95 | 68 |
| 384 | 121 | 70 | 49 |
| 768 | 161 | 56 | 44 |

*(II) Transcritical flow*

The steady-state surface level $h + b$ for the transcritical flow computed by NMGM-3 on a 96-cell grid for the same channel configurations, together with the reference analytical solution from [27] and the bottom topography, is plotted in Figure 7. Compared with the subcritical case, the flow undergoes a critical transition near the bump, resulting in a more nonlinear free-surface profile. Nevertheless, the numerical solution remains in close agreement with the reference for all channel configurations, indicating that NMGM-3 retains its accuracy for the more challenging transcritical regime.

In Table 4, the number of Newton iterations required by NMGM-3 to reach the steady state on successively refined grids is presented. As in the subcritical case, the solver converges robustly for all tested transcritical flows, although generally more Newton iterations are required, which reflects the increased difficulty of the problem due to the critical transition near the bump. For the channel with constant width, the number of Newton iterations remains relatively stable under grid refinement. A similar grid-independent behavior is observed for the channel with a left-shifted contraction. By contrast, for the channel with a right-shifted contraction, the number of Newton iterations increases rapidly as the grid is refined, yet convergence is still achieved.

It is worth noting that the rapid increase in iteration numbers in the right-shifted contraction case appears to be partly related to the perturbation parameter $\epsilon$ used in the numerical differentiation (12) for Jacobian evaluation. In practice, this increase can be mitigated by adjusting $\epsilon$. For instance, on the 768-cell grid, reducing $\epsilon$ from 0.05 to 0.02 noticeably decreases the iteration numbers from 537 to 269.

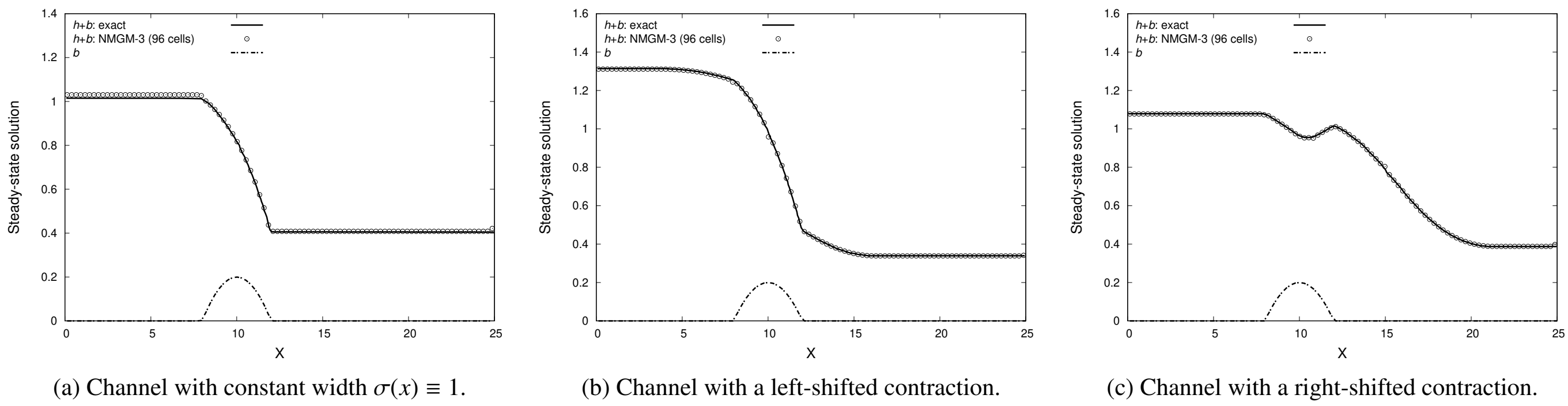

(a) Channel with constant width $\sigma(x) \equiv 1$.
(b) Channel with a left-shifted contraction.
(c) Channel with a right-shifted contraction.

Figure 7: Transcritical flow over a bump: steady-state surface level $h + b$ obtained by NMGM-3 on a 96-cell grid, compared with the exact solution from [27].

Table 4: Transcritical flow over a bump: number of Newton iterations required by NMGM-3 for different channel widths on successively refined grids.

| No. of cells | Constant width | Left-shifted contraction | Right-shifted contraction |
| --- | --- | --- | --- |
| 96 | 144 | 209 | 208 |
| 192 | 166 | 215 | 256 |
| 384 | 145 | 245 | 370 |
| 768 | 161 | 200 | 537 |

### 5.2. *Two-dimensional case*

Simulations for the two-dimensional problems are carried out on three successively refined grids with resolutions of $16 \times 8$, $32 \times 16$, and $64 \times 32$ cells. In these cases, characteristic decomposition [19] is employed within the WENO reconstruction to reduce spurious oscillations, and the perturbation parameter $\epsilon$ in the numerical differentiation (12) is fixed at $10^{-6}$. To maintain robust multigrid convergence across different grid levels, the numbers of pre- and post-smoothing steps are set to $\nu_1 = \nu_2 = 10$, 15, and 20 for the three grids, respectively.

#### 5.2.1. *Subcritical flow over a parabolic bump*

In this example, we consider a two-dimensional subcritical flow over a parabolic bump in the rectangular domain $[-4, 4] \times [0, 4]$. The bottom topography is defined as

$$b(x, y) = \begin{cases} 0.4 - 0.4\left(x^2 + (y-2)^2\right), & x^2 + (y-2)^2 < 1, \\ 0, & \text{otherwise}. \end{cases} \tag{29}$$

At the inflow boundary $x = -4$, the discharge $hu = 1$ is prescribed, whereas at the outflow boundary $x = 4$, the water depth $h = 1$ is imposed. At the sidewalls $y = 0$ and $y = 4$, reflective boundary conditions are applied to represent the impermeable slip walls.

The NMGM solver is initialized on the coarsest auxiliary grid of $8 \times 4$ cells with the uniform state

$$h(x, y) = 1, \quad hu(x, y) = 1, \quad hv(x, y) = 0. \tag{30}$$

As in the one-dimensional treatment, the initial guess on successively refined grids is obtained by interpolating the steady-state solution from the preceding coarser grid.

Figure 8 shows a contour plot of the steady-state free surface $h + b$ computed by NMGM-9 on a $64 \times 32$ grid. The numerical solution captures the overall flow structure over the parabolic bottom bump without introducing spurious oscillations, demonstrating the stability and non-oscillatory behavior of the solver.

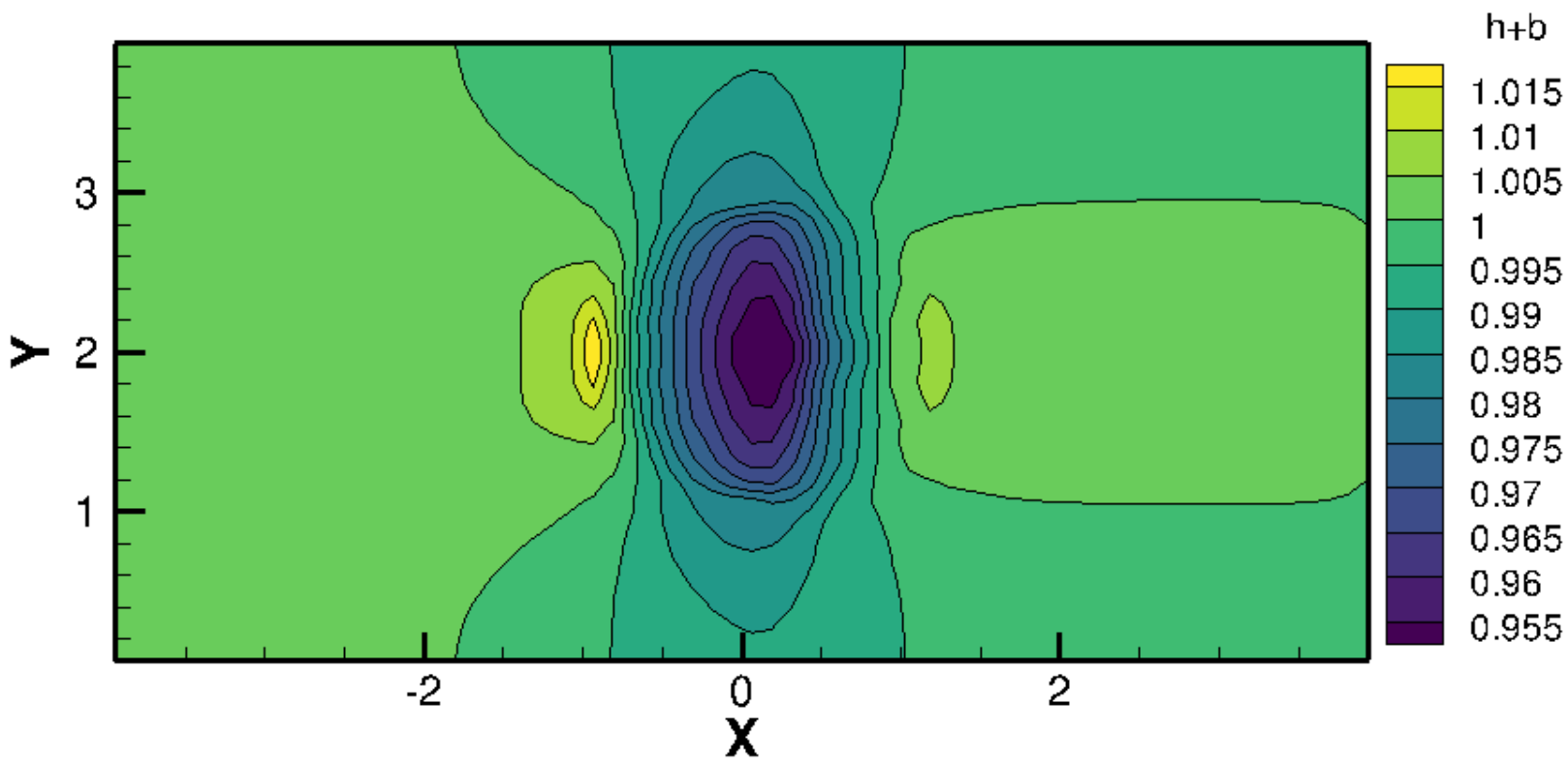

Figure 8: Two-dimensional subcritical flow over a parabolic bump: contour plot of the steady-state surface $h + b$ obtained by NMGM-9 on a $64 \times 32$-cell grid.

The convergence behavior of NMGM-9 on the three tested grids is shown in Figure 9, where the $\ell^1$-norm of the global residual is plotted against the number of Newton iterations. In all cases, the residual decreases steadily, indicating robust convergence of the proposed solver. As the grid is refined, the convergence rate becomes slightly slower. This leads to an increase in the total number of Newton iterations. Nevertheless, the growth in iteration numbers remains moderate under grid refinement, suggesting favorable scalability of the solver.

Simulations on the three grids are also carried out using NMGM-21 with the same parameters as those employed for NMGM-9. Here, NMGM-9 and NMGM-21 correspond to the simplified and full Jacobian formulations, respectively, in analogy with the one-dimensional solvers NMGM-3 and NMGM-5. Consistent with the one-dimensional results, the two solvers show nearly identical convergence behavior, requiring almost the same number of Newton iterations to reach the steady state in all cases. Accordingly, their main difference in computational cost lies in the construction of the Jacobian matrix. To assess this efficiency difference, Table 5 reports the CPU time required for assembling the Jacobian once on the three grids. The results indicate that NMGM-9 consistently incurs a substantially lower cost than NMGM-21, yielding a speedup close to a factor of two, which is in agreement with the theoretical reduction in the Jacobian construction cost predicted in subsection 4.1, thereby confirming the improved efficiency of the simplified Jacobian approach.

Table 5: Two-dimensional subcritical flow over a parabolic hump: CPU time (s) for Jacobian construction using NMGM-21 and NMGM-9.

| Solver | No. of cells | | |
|---|---|---|---|
| | $16 \times 8$ | $32 \times 16$ | $64 \times 32$ |
| NMGM-21 | 11 | 53 | 293 |
| NMGM-9 | 5 | 25 | 144 |

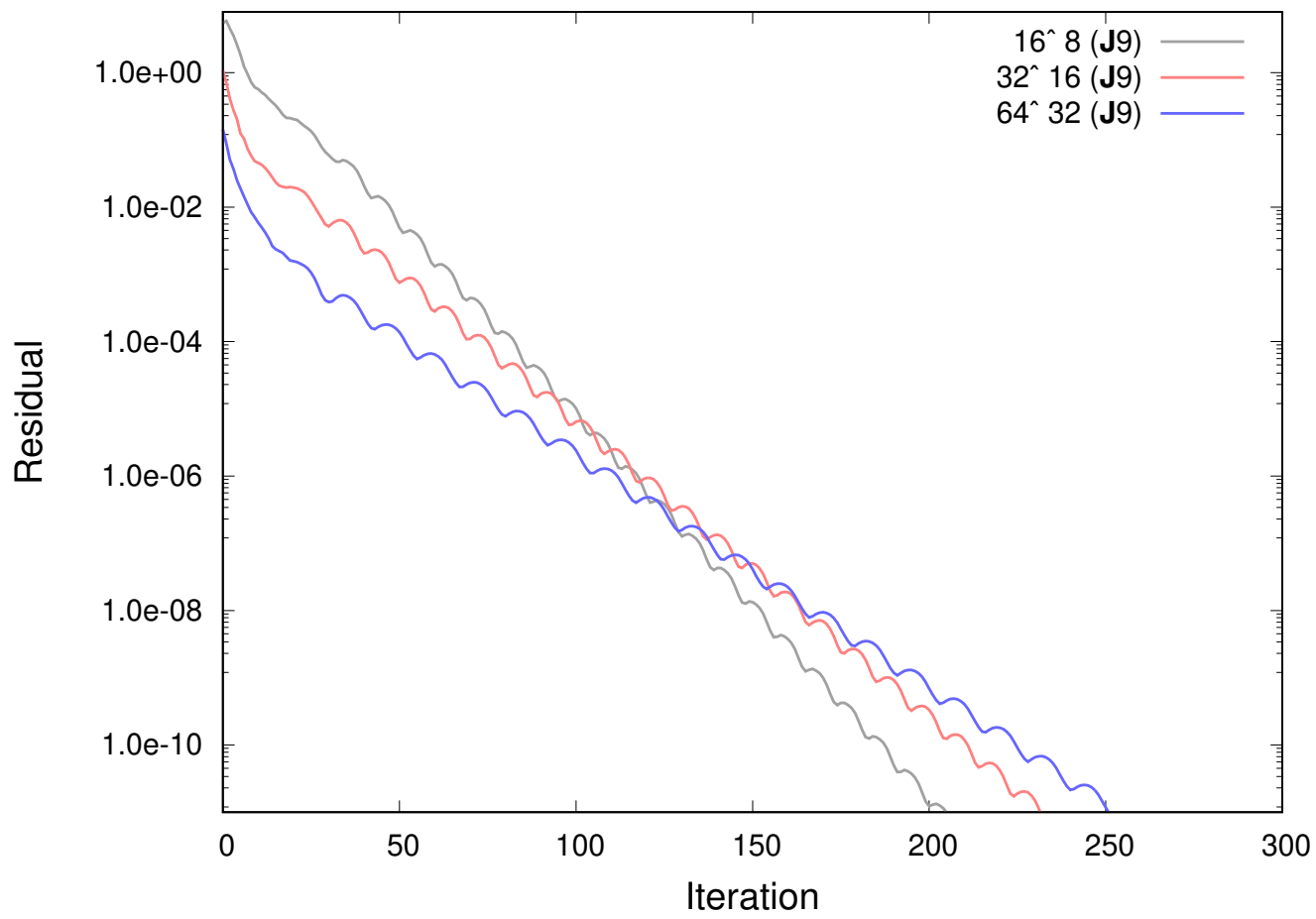

Figure 9: Two-dimensional subcritical flow over a parabolic hump: convergence histories of NMGM-9 on three tested grids.

### 5.2.2. *Hydraulic jump over a wedge problem*

The last example considers a supercritical flow over a wedge with an angle of $\theta_w = 8.95°$, which induces a hydraulic jump within the flow domain. This problem has been investigated in [5, 28, 29, 30] to assess the conservative and non-oscillatory properties of numerical schemes for shallow water flows with strong discontinuities.

Consistent with the setup in [30], the computational domain is $[0, 5]\times[0, 2]$. On the left and top boundaries, inflow conditions are prescribed with the water depth $h = 1$ and velocity $(u, v) = 8.57(\cos(-\theta_w), \sin(-\theta_w))$. On the right boundary $x = 5$, an outflow condition is applied by extrapolating the conservative variables from the interior. On the bottom boundary $y = 0$, the outflow condition is applied for $x < 1/\cos(\theta_w)$, whereas a reflective boundary condition is imposed on the remaining portion to represent the solid wedge surface.

For the NMGM solver, the initial guess is constructed on the coarsest grid of $16 \times 8$ cells using uniform states

$$h^{(0)}(x, y) = 1, \quad u^{(0)}(x, y) = 8.57\cos(-\theta_w), \quad v^{(0)}(x, y) = 8.57\sin(-\theta_w). \tag{31}$$

On finer grids, the initial solution is obtained by interpolation from the converged solution on the preceding coarser grid.

As shown in Figure 10, the steady-state water depth $h$ computed by NMGM-9 on a $64 \times 32$ grid exhibits a clear oblique hydraulic jump inside the computational domain, which is consistent with the results reported in [30]. Specifically, a shallow, nearly uniform region occupies the upper-left part of the domain, while a deeper region appears in the lower-right part. These regions correspond to the supercritical and subcritical flow states, respectively. The transition between them is aligned along an approximately straight line, and the jump is captured within a relatively narrow zone without noticeable spurious oscillations, illustrating the non-oscillatory behavior of the NMGM solver in the presence of a strong discontinuity.

The convergence histories of NMGM-9 as well as NMGM-21 on the three tested grids are shown in Figure 11, where the $\ell^1$-norm of the global residual is plotted against the number of Newton iterations. As observed in the previous example, the residual decreases steadily in all cases, indicating robust convergence despite the presence of a strong discontinuity. It can also be seen that the NMGM solver requires more Newton iterations on the coarsest grid than on the finer grids. This is mainly because the initial guess on the coarsest grid is constructed directly using (31), rather than obtained by interpolation from a coarser-grid solution, leading to a slower residual reduction during the initial iterations. Additionally, NMGM-9 and NMGM-21 exhibit nearly identical convergence histories on all three grids. This indicates that the simplified Jacobian $\boldsymbol{J}9$ preserves the essential nonlinear coupling required by NMGM even in the presence of a strong hydraulic jump. Consequently, comparable convergence can be achieved with the cheaper Jacobian approximation, implying improved computational efficiency without compromising robustness.

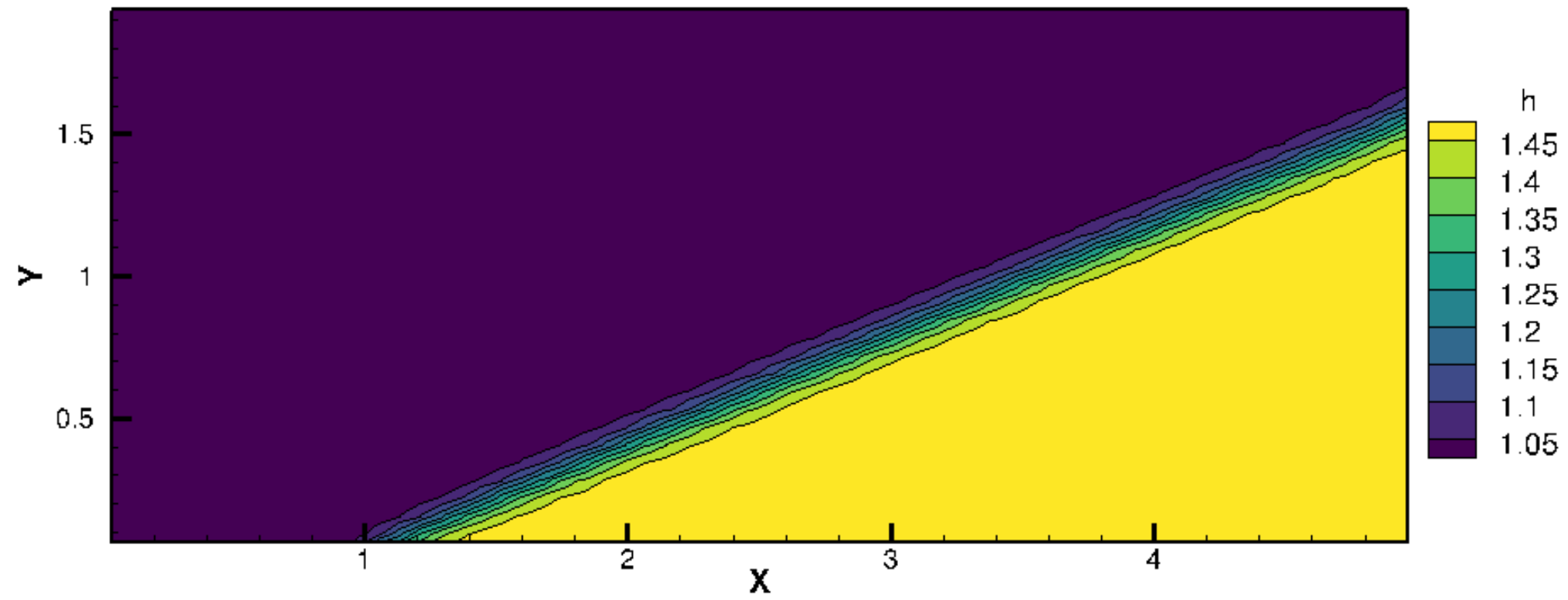

Figure 10: Hydraulic jump over a wedge problem: contour plot of the steady-state water depth $h$ obtained by NMGM-9 on a $64 \times 32$-cell grid.

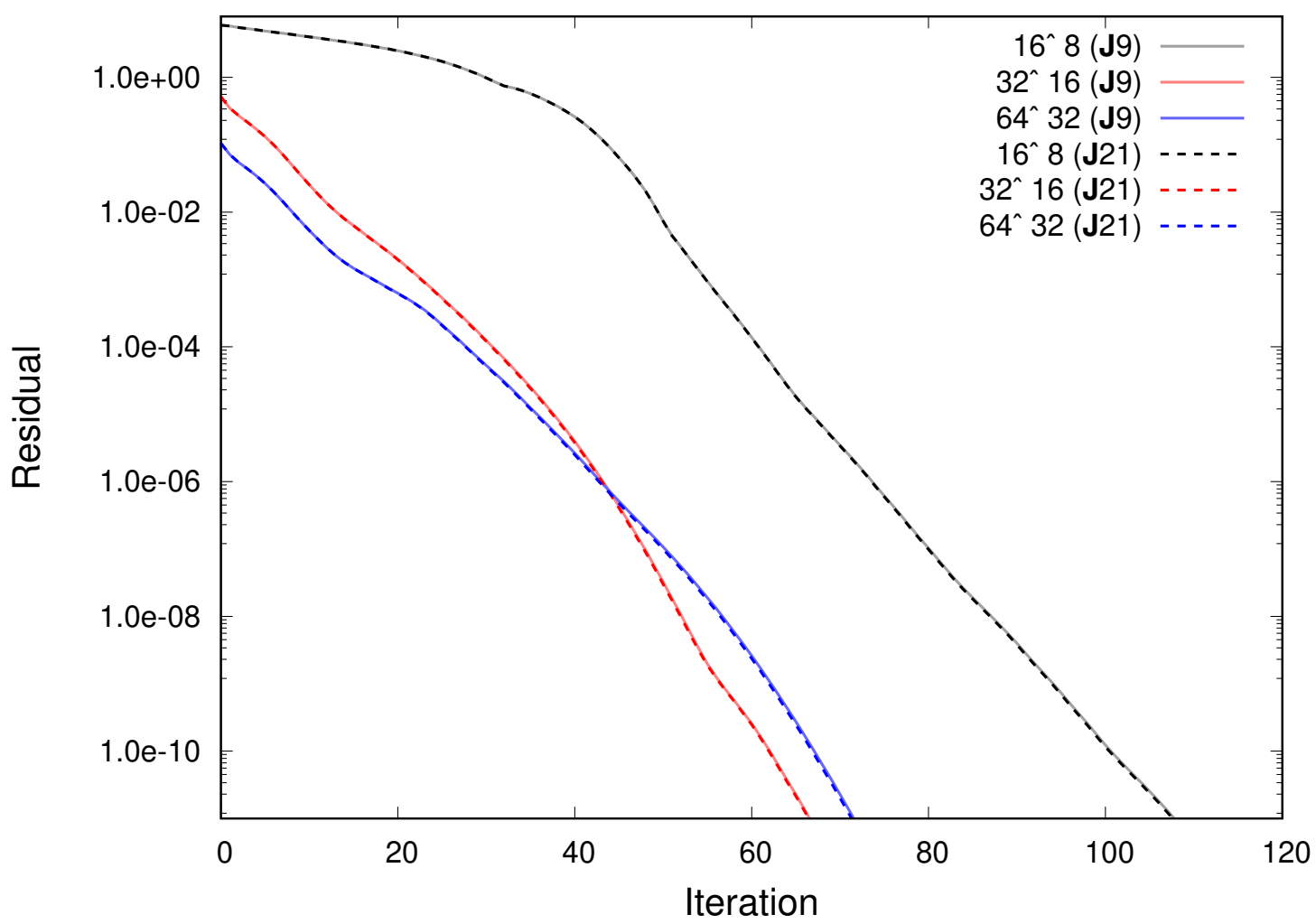

Figure 11: Hydraulic jump over a wedge problem: convergence histories of NMGM-9 and NMGM-21 on three tested grids.

## 6. Conclusion

In this paper, a third-order Newton multigrid method (NMGM) has been developed for the steady-state shallow water equations in open channels. The method combines a high-order finite volume WENO discretization with an efficient Newton multigrid framework for solving the resulting nonlinear system. To alleviate the high computational cost associated with the Jacobian matrix in high-order schemes, a simplified Jacobian approximation based on reduced stencils has been proposed for both one- and two-dimensional problems. This approach significantly reduces the cost of Jacobian construction while retaining the essential nonlinear coupling required for robust and fast convergence. A number of numerical experiments have been conducted to demonstrate the accuracy, efficiency and robustness of the proposed method. The results show that the simplified Jacobian approach achieves substantial computational savings without compromising accuracy or robustness across all tests. Specifically, it reduces the computational cost by more than 30% in one dimension and 50% in two dimension, while exhibiting nearly identical convergence histories, in comparison to the full Jacobian approach. Overall, the proposed NMGM provides an efficient and robust framework for steady-state shallow water simulations and shows considerable potential for extension to more complex flow problems.

## Acknowledgements

This work was partially supported by the National Natural Science Foundation of China, No. 12171240, the National Science and Technology Major Project, No. J2019-II-0007-0027, and the Postgraduate Research & Practice Innovation Program of Jiangsu Province, No. KYCX24_0523. The computational resources were supported by High Performance Computing Platform of Nanjing University of Aeronautics and Astronautics, China.

## Appendix A. Computation of the initial cell averages from the coarser-level solution

For convenience, we denote the variables and operators related to the coarse and fine grid levels by subscripts $L$ and $l$, respectively. Given the coarse-level solution $\overline{\boldsymbol{U}}_L$, our goal is to compute the fine-level cell averages $\overline{\boldsymbol{U}}_l$ in a conservative and high-order accurate manner.

In the one-dimensional case, the quadratic polynomial that interpolates the neighboring coarse-level cell averages $\overline{U}_{j-1,L}$, $\overline{U}_{j,L}$, and $\overline{U}_{j+1,L}$ reads

$$P_{j,L}(x) = \overline{U}_{j,L} + \frac{\overline{U}_{j+1,L} - \overline{U}_{j-1,L}}{2}\left(\frac{x - x_{j,L}}{\Delta x_L}\right) + \frac{\overline{U}_{j-1,L} - 2\overline{U}_{j,L} + \overline{U}_{j+1,L}}{2}\left[\left(\frac{x - x_{j,L}}{\Delta x_L}\right)^2 - \frac{1}{12}\right]. \tag{A.1}$$

Note that $I_{j,L} = I_{2j-1,l} \bigcup I_{2j,l}$, we obtain the fine-level cell averages over $I_{2j-1,l}$ and $I_{2j,l}$ by

$$\begin{aligned} \overline{U}_{2j-1,l} &= \frac{1}{\Delta x_l}\int_{I_{2j-1,l}} P_{j,L}(x)\,\mathrm{d}x = \frac{\overline{U}_{j-1,L}}{8} + \overline{U}_{j,L} - \frac{\overline{U}_{j+1,L}}{8}, \\ \overline{U}_{2j,l} &= \frac{1}{\Delta x_l}\int_{I_{2j,l}} P_{j,L}(x)\,\mathrm{d}x = -\frac{\overline{U}_{j-1,L}}{8} + \overline{U}_{j,L} + \frac{\overline{U}_{j+1,L}}{8}. \end{aligned} \tag{A.2}$$

In the two-dimensional case, the interpolating polynomial $P_{ij,L}(x,y)$ associated with the $(i,j)$th coarse-level cell $I_{ij,L}$ is obtained through a dimension-by-dimension reconstruction. Specifically, we have

$$P_{ij,L}(x,y) = P_0(x) + \frac{P_1(x) - P_{-1}(x)}{2}\left(\frac{y - y_{j,L}}{\Delta y_L}\right) + \frac{P_{-1}(x) - 2P_0(x) + P_1(x)}{2}\left[\left(\frac{y - y_{j,L}}{\Delta y_L}\right)^2 - \frac{1}{12}\right], \tag{A.3}$$

where $P_k(x)$, $k = -1, 0, 1$, are quadratic polynomials interpolating the coarse-level cell averages $\overline{U}_{i-1,j+k,L}$, $\overline{U}_{i,j+k,L}$, and $\overline{U}_{i+1,j+k,L}$ as defined by the one-dimensional formula (A.1). Then the fine-level cell averages over $I_{2i-1,2j-1,l}$, $I_{2i-1,2j,l}$, $I_{2i,2j-1,l}$ and $I_{2i,2j,l}$ are obtained by integrating $P_{ij,L}(x,y)$ over each cell as

$$\begin{aligned} \overline{U}_{2i-1,2j-1,l} &= \frac{1}{\Delta x_l \Delta y_l}\int_{I_{2i-1,2j-1,l}} P_{ij,L}(x,y)\,\mathrm{d}x\,\mathrm{d}y = \sum_{k=-1}^{1} \frac{8(1+k)(1-k) - k}{64}\left(\overline{U}_{i-1,j+k,L} + 8\overline{U}_{i,j+k,L} - \overline{U}_{i+1,j+k,L}\right), \\ \overline{U}_{2i-1,2j,l} &= \frac{1}{\Delta x_l \Delta y_l}\int_{I_{2i-1,2j,l}} P_{ij,L}(x,y)\,\mathrm{d}x\,\mathrm{d}y = \sum_{k=-1}^{1} \frac{8(1+k)(1-k) + k}{64}\left(\overline{U}_{i-1,j+k,L} + 8\overline{U}_{i,j+k,L} - \overline{U}_{i+1,j+k,L}\right), \\ \overline{U}_{2i,2j-1,l} &= \frac{1}{\Delta x_l \Delta y_l}\int_{I_{2i,2j-1,l}} P_{ij,L}(x,y)\,\mathrm{d}x\,\mathrm{d}y = \sum_{k=-1}^{1} \frac{8(k+1)(k-1) + k}{64}\left(\overline{U}_{i-1,j+k,L} - 8\overline{U}_{i,j+k,L} - \overline{U}_{i+1,j+k,L}\right), \\ \overline{U}_{2i,2j,l} &= \frac{1}{\Delta x_l \Delta y_l}\int_{I_{2i,2j,l}} P_{ij,L}(x,y)\,\mathrm{d}x\,\mathrm{d}y = \sum_{k=-1}^{1} \frac{8(k+1)(k-1) - k}{64}\left(\overline{U}_{i-1,j+k,L} - 8\overline{U}_{i,j+k,L} - \overline{U}_{i+1,j+k,L}\right). \end{aligned} \tag{A.4}$$